\documentclass[11pt]{amsart}
\usepackage{amssymb,amsfonts,amsthm,amscd,amsmath,amsgen,upref}

\usepackage[margin=1in]{geometry}

\theoremstyle{plain}
\newtheorem{cor}{Corollary}

\newtheorem{prop}[cor]{Proposition}
\newtheorem{con}[cor]{Control}
\newtheorem{thm}[cor]{Theorem}
\theoremstyle{definition}

\numberwithin{cor}{section}
\numberwithin{equation}{section}

\DeclareMathOperator{\tr}{tr}

\DeclareMathOperator{\C}{C}
\DeclareMathOperator{\USC}{USC}
\DeclareMathOperator{\LSC}{LSC}
\DeclareMathOperator{\BUC}{BUC}
\DeclareMathOperator{\Lip}{Lip}

\DeclareMathOperator{\Supp}{Supp}

\renewcommand{\d}{d} 

\newcommand{\abs}[1]{\lvert#1\rvert}
\newcommand{\norm}[1]{\lVert#1\rVert}

\def\XXint#1#2#3{{\setbox0=\hbox{$#1{#2#3}{\int}$ }
\vcenter{\hbox{$#2#3$ }}\kern-.6\wd0}}

\title{Exit Laws of Isotropic Diffusions in Random Environment from Large Domains}

\author{Benjamin J. Fehrman$^\dagger$}
\thanks{$^\dagger$ This material is based upon work supported by the National Science Foundation Mathematical Sciences Postdoctoral Research Fellowship under Grant Number 1502731.}

\address{Max Planck Institute for Mathematics in the Sciences \\ Inselstra{$\ss$}e 22\\ 04103 Leipzig, Germany}
\email{fehrman@mis.mpg.de}

\date{23/1/2016}

\begin{document}

\maketitle

\begin{abstract}  This paper studies, in dimensions greater than two, stationary diffusion processes in random environment which are small, isotropic perturbations of Brownian motion satisfying a finite range dependence.  Such processes were first considered in the continuous setting by Sznitman and Zeitouni \cite{SZ}.  Building upon their work, it is shown by analyzing the associated elliptic boundary-value problem that, almost surely, the smoothed (in the sense that the boundary data is continuous) exit law of the diffusion from large domains converges, as the domain's scale approaches infinity, to that of a Brownian motion.  Furthermore, a rate for the convergence is established in terms of the modulus of the boundary condition.\end{abstract}

\section{Introduction}\label{introduction}

The purpose of this paper is to characterize the smoothed exit distributions from large domains associated to the diffusion in random environment determined by the generator \begin{equation}\label{intro_generator}\frac{1}{2}\sum_{i,j=1}^da_{ij}(x,\omega)\frac{\partial^2}{\partial x_i\partial x_j}+\sum_{i=1}^db_i(x,\omega)\frac{\partial}{\partial x_i},\end{equation} where the environment, as described by a uniformly elliptic diffusion matrix $A(x,\omega)=(a_{ij}(x,\omega))$ and drift $b(x,\omega)=(b_i(x,\omega))$, is indexed by an underlying probability space $(\Omega,\mathcal{F},\mathbb{P})$.

The family of stochastic processes associated to the generator (\ref{intro_generator}) will be assumed to be a stationary, isotropic perturbation of Brownian motion satisfying a finite range dependence.  Specifically, there exists a measure-preserving group of transformations $\left\{\tau_x\right\}_{x\in\mathbb{R}^d}$ such that, for each $x,y\in\mathbb{R}^d$ and $\omega\in\Omega$, $$A(x+y,\omega)=A(x,\tau_y\omega)\;\;\textrm{and}\;\;b(x+y,\omega)=b(x,\tau_y\omega).$$  Whenever subsets $A,B\subset\mathbb{R}^d$ are sufficiently distant in space, the sigma algebras $$\sigma(A(x,\cdot), b(x,\cdot)\;|\;x\in A)\;\;\textrm{and}\;\;\sigma(A(x,\cdot),b(x,\cdot)\;|\;x\in B)\;\;\textrm{are independent.}$$  For every orthogonal transformation $r$ of $\mathbb{R}^d$ preserving the coordinate axis, for each $x\in\mathbb{R}^d$, the random variables $$(A(rx,\omega),b(rx,\omega))\;\;\textrm{and}\;\;(rA(x,\omega)r^t,rb(x,\omega))\;\;\textrm{have the same law.}$$  Finally, for a constant $\eta>0$ to be chosen small, for every $x\in\mathbb{R}^d$ and $\omega\in\Omega$, $$\abs{A(x,\omega)-I}<\eta\;\;\textrm{and}\;\;\abs{b(x,\omega)}<\eta.$$  Such environments were considered in the continuous setting by Sznitman and Zeitouni \cite{SZ} and correspond to the analogue of the discrete framework studied by Bricmont and Kupiainen \cite{BK}.

The exit distribution of these processes will be understood by analyzing the asymptotic behavior, as $\epsilon\rightarrow 0$, of solutions \begin{equation}\label{intro_eq}\left\{\begin{array}{ll} \frac{1}{2}\tr(A(x,\omega)D^2v^\epsilon)+b(x,\omega)\cdot Dv^\epsilon=0 & \textrm{on}\;\;U/\epsilon, \\ v^\epsilon=f(\epsilon x) & \textrm{on}\;\;\partial U/\epsilon,\end{array}\right.\end{equation}  where, in dimension at least two, the domain $U\subset\mathbb{R}^d$ satisfies an exterior ball condition and $f\in C(\partial U)$ is continuous on the boundary.

Writing $E_{x,\omega}$ for the expectation describing the diffusion $X_t$ associated to the generator (\ref{intro_generator}) beginning from $x$ in environment $\omega$, and $\tau^\epsilon$ for the exit time from $U/\epsilon$, solutions of (\ref{intro_eq}) admit the representation \begin{equation}\label{intro_rep} v^\epsilon(x)=E_{x,\omega}(f(\epsilon X_{\tau^{\epsilon}}))\;\;\textrm{on}\;\;\overline{U}/\epsilon.\end{equation}  The continuity of the boundary data $f$ in this setting corresponds to a necessary smoothing of the exit distribution since, and as described below, the presence of traps (which, loosely speaking, mean portions of space where the drift has a strong effect) along the boundary preclude, in the case of discontinuous boundary data, an almost sure characterization of the exit measures defining the solutions $v^\epsilon$ in the limit.

Observe, following a rescaling, that $u^\epsilon(x)=v^\epsilon(x/\epsilon)$ satisfies \begin{equation}\label{intro_res}\left\{\begin{array}{ll}\frac{1}{2}\tr(A(\frac{x}{\epsilon},\omega)D^2u^\epsilon)+\frac{1}{\epsilon}b(\frac{x}{\epsilon},\omega)\cdot Du^\epsilon=0 & \textrm{on}\;\;U, \\ u^\epsilon=f(x) & \textrm{on}\;\;\partial U,\end{array}\right.\end{equation} and, in view of (\ref{intro_rep}), $$u^\epsilon(x)=E_{\frac{x}{\epsilon},\omega}(f(\epsilon X_{\tau_\epsilon}))=E_{\frac{x}{\epsilon},\omega}(f(\epsilon X_{\frac{\epsilon^2\tau_\epsilon}{\epsilon^2}})),$$ where $\epsilon^2\tau_\epsilon$ is by definition the exit time from $U$ of the rescaled process $\epsilon X_{t/\epsilon^2}$.  The behavior of this rescaled process on $\mathbb{R}^d$ was characterized in \cite{SZ}, where it was proven that, on a subset of full probability and for a deterministic $\overline{\alpha}>0$, \begin{equation}\label{intro_sz}\epsilon X_{\cdot/\epsilon^2}\;\textrm{converges as}\;\epsilon\rightarrow 0\;\textrm{in law on}\;\mathbb{R}^d\;\textrm{to a Brownian motion with variance}\;\overline{\alpha}.\end{equation}  See Section \ref{preliminaries} for the precise statement and details.

The primary aim of this paper is to establish the analogous result for the exit distribution, and this is achieved by characterizing, on a subset of full probability, the asymptotic behavior as $\epsilon\rightarrow 0$ of solutions to (\ref{intro_res}).

\begin{thm}\label{intro_main} There exists a subset of full probability such that, for every bounded domain $U\subset\mathbb{R}^d$ satisfying an exterior ball condition, the solutions of (\ref{intro_eq}) converge uniformly on $\overline{U}$, as $\epsilon\rightarrow 0$, to the solution \begin{equation}\label{intro_hom}\left\{\begin{array}{ll} \Delta \overline{u}=0 & \textrm{on}\;\;U, \\ \overline{u}=f(x) & \textrm{on}\;\;\partial U.\end{array}\right.\end{equation}\end{thm}

The proof relies strongly upon the results obtained in \cite{SZ}, and in particular between a comparison there obtained, with scaling analogous to (\ref{intro_eq}), between solutions of the parabolic equation \begin{equation}\label{intro_parabolic} \left\{\begin{array}{ll} u^\epsilon_t=\frac{1}{2}\tr(A(x,\omega)D^2u^\epsilon)+b(x,\omega)\cdot Du^\epsilon & \textrm{on}\;\;\mathbb{R}^d\times(0,\infty), \\ u^\epsilon=f(\epsilon x) & \textrm{on}\;\;\mathbb{R}^d\times\left\{0\right\}, \end{array}\right.\end{equation} and the parabolic analogue of (\ref{intro_hom}) with high probability and on large scales in space and time.  The details of this argument are presented in Section \ref{inductive}.

The comparison is used, as it was in \cite{SZ}, to construct a coupling between the process defined by the generator (\ref{intro_generator}) and a Brownian motion with variance approximately $\overline{\alpha}$ along a discrete sequence of time steps.  See Proposition \ref{couple_main} in Section \ref{coupling}.  And, this coupling allows for the introduction of a discrete version of equation (\ref{intro_eq}).  Namely, the process is evaluated along the aforementioned discrete sequence of time steps and stopped as soon as it hits a neighborhood of the compliment of the domain.

However, the approximation suggested here is typically insufficient to characterize the limiting behavior of solutions to (\ref{intro_eq}) and its rescaling (\ref{intro_res}) since the time steps are not sufficiently fine to preclude the emergence of traps created by the drift.  And, in this setting, the traps are twofold.  Considering the process associated to (\ref{intro_res}), and corresponding to the generator \begin{equation}\label{intro_res_gen}\frac{1}{2}\sum_{i,j=1}^da_{ij}(\frac{x}{\epsilon},\omega)\frac{\partial^2}{\partial x_i \partial x_j}+\frac{1}{\epsilon}\sum_{i=1}^db_i(\frac{x}{\epsilon},\omega)\frac{\partial}{\partial x_i},\end{equation} the presence of the singular (in $\frac{1}{\epsilon}$) drift can first act to confine the particle to create, in expectation, an exponentially growing (in $\frac{1}{\epsilon}$) exit time.  The probability that the exit time is large is first controlled, though sub-optimally, by Proposition \ref{exit_main} in Section \ref{section_exit_time}.

Second, the drift can repel the process from the boundary, and thereby make impossible the existence (in general) of barriers which are effective at scales greater than $\epsilon$.  This difficulty is overcome by combining, in Section \ref{main}, the coupling obtained in Proposition \ref{couple_main} of Section \ref{coupling} with estimates concerning the exit time of Brownian motion from the slightly inflated domains $$U_\delta=\left\{\;x\in\mathbb{R}^d\;|\;d(x,U)<\delta\;\right\},$$ where $\delta\rightarrow 0$ as $\epsilon\rightarrow 0$.  These estimates are proven in Propositions \ref{disc_u} and \ref{disc_u_scale} of Section \ref{Brownian_exit}.  Then, at points near the boundary $\partial U$, the exit of the Brownian motion from a somewhat larger domain of the type $U_\delta$ is shown to compel, with high probability, the exit of the diffusion in random environment from $U$.  See Proposition \ref{end_time} of Section \ref{main}.  It is this fact that establishes the efficacy of the discrete approximation and ultimately the proof of Theorem \ref{intro_main}.

Finally, in Section \ref{section_rate}, the convergence established in Theorem \ref{intro_main} is made quantitative assuming first that the boundary data $f$ is the restriction of a bounded, uniformly continuous function on $\mathbb{R}^d$.  Namely, \begin{equation}\label{intro_rate_boundary}\textrm{assume}\;f\in\BUC(\mathbb{R}^d)\;\;\textrm{with modulus}\;\sigma_f\;\textrm{satisfying}\;\abs{f(x)-f(y)}\leq\sigma_f(\abs{x-y})\end{equation} for all $x,y\in\mathbb{R}^d$.  

\begin{thm}\label{intro_quantitative} Assume (\ref{intro_rate_boundary}).  There exist constants $0<c_0,c_1<1$ and $C>0$ such that, on a subset of full probability, for all $\epsilon>0$ sufficiently small depending on $\omega$, the solutions of (\ref{intro_res}) and (\ref{intro_hom}) satisfy $$\norm{u^\epsilon-\overline{u}}_{L^\infty(\overline{U})}\leq C\norm{f}_{L^\infty(\mathbb{R}^d)}\epsilon^{c_0}+C\sigma_f(\epsilon^{c_1}).$$\end{thm}

A standard extension argument then allows Theorem \ref{intro_quantitative} to be extended to arbitrary continuous functions on the boundary, provided the domain is smooth.  In this case, \begin{equation}\label{intro_rate_boundary_2} \textrm{assume the domain}\;U\;\textrm{is smooth},\end{equation} and \begin{equation}\label{intro_rate_boundary_1}\textrm{assume}\;f\in\C(\partial U)\;\;\textrm{with modulus}\;\sigma_f\;\textrm{satisfying}\;\abs{f(x)-f(y)}\leq\sigma_f(\abs{x-y})\end{equation} for all $x,y\in\partial U$.

\begin{thm}\label{intro_quantitative_1} Assume (\ref{intro_rate_boundary_2}) and (\ref{intro_rate_boundary_1}).  There exist constants $0<c_0,c_1<1$, $C_1=C_1(U)>0$ depending upon the domain and $C>0$ such that, on a subset of full probability, for all $\epsilon>0$ sufficiently small depending on $\omega$, the solutions of (\ref{intro_res}) and (\ref{intro_hom}) satisfy $$\norm{u^\epsilon-\overline{u}}_{L^\infty(\overline{U})}\leq C\norm{f}_{L^\infty(\partial U)}\epsilon^{c_0}+C\sigma_f(C_1\epsilon^{c_1}).$$\end{thm}

Diffusion processes on $\mathbb{R}^d$ in the stationary ergodic setting were first considered in the case $b(x,\omega)=0$ by Papanicolaou and Varadhan \cite{PV1}.  Furthermore, in the case that (\ref{intro_eq}) can be rewritten in divergence form, these diffusions and associated boundary value problems were studied in Papanicolaou and Varadhan \cite{PV}, and further results have been obtained by De Masi, Ferrari, Goldstein and Wick \cite{MFGW}, Kozlov \cite{Kozlov}, Olla \cite{Olla} and Osada \cite{Osada}.  However, for general drifts $b(x,\omega)$ which are neither divergence free nor a gradient of a stationary field, considerably less is known.

Indeed, the results of \cite{SZ}, which apply to the isotropic, perturbative regime described above, and later extended by Fehrman \cite{F2, F1} are the only such available.  And, to this point, the characterization of the asymptotic behavior of boundary value problems like (\ref{intro_res}) has remained open.  However, some results do exist for the analogous discrete framework.  Bolthausen and Zeitouni \cite{BZ} characterized the exit distributions from large balls (so, taking $U=B_1$) of random walks in random environment which are small, isotropic perturbations of a simple random walk, and their work was later refined by Baur and Bolthausen \cite{BB} under a somewhat less stringent isotropy assumption.  Finally, Baur \cite{Baur} has recently obtained results concerning the exit time from large balls of processes satisfying a quenched symmetry assumption along a single coordinate direction.

The methods of this paper differ significantly from those of \cite{BB,Baur,BZ}, which develop an induction scheme to propagate estimates concerning the convergence of the exit law of the diffusion in random environment to the uniform measure on the boundary of the ball, by instead relying upon the results of \cite{SZ} obtained in the parabolic setting.  Furthermore, these methods apply to arbitrary bounded domains satisfying an exterior ball condition.

The paper is organized so that, in Section \ref{preliminaries}, the notation and assumptions are presented and, in Section \ref{inductive}, the most relevant aspects of \cite{SZ} are reviewed and the primary probabilistic statement concerning the random environment is presented.  In Section \ref{section_exit_time}, the exit time of the process in random environment is controlled in probability, and the global coupling between the process in random environment and Brownian motion is constructed in Section \ref{coupling}.  The exit time of Brownian motion at points near the boundary of the inflated domains $U_\delta$ is controlled in Section \ref{Brownian_exit}, and the efficacy of the discrete approximation, as defined through the coupling, and ultimately the proof of Theorem \ref{intro_main} are presented in Section \ref{main}.  Finally, the rates of convergence in Theorems \ref{intro_quantitative} and \ref{intro_quantitative_1} appear in Section \ref{section_rate}.

\subsection*{Acknowledgments}

I would like to thank Professors Panagiotis Souganidis and Ofer Zeitouni for many useful conversations.

\section{Preliminaries}\label{preliminaries}

\subsection{Notation}

Elements of $\mathbb{R}^d$ and $[0,\infty)$ are denoted by $x$ and $y$ and $t$ respectively and $(x,y)$ denotes the standard inner product on $\mathbb{R}^d$.  The gradient in space and derivative in time of a scalar function $v$ are written $Dv$ and $v_t$, while $D^2v$ stands for the Hessian of $v$.  The spaces of $k\times l$ and $k\times k$ symmetric matrices with real entries are respectively written $\mathcal{M}^{k\times l}$ and $\mathcal{S}(k)$.  If $M\in\mathcal{M}^{k\times l}$, then $M^t$ is its transpose and $\abs{M}$ is its norm $\abs{M}=\tr(MM^t)^{1/2}.$  If $M$ is a square matrix, the trace of $M$ is written $\tr(M)$.  The Euclidean distance between subsets $A,B\subset\mathbb{R}^d$ is $$d(A,B)=\inf\left\{\;\abs{a-b}\;|\;a\in A, b\in B\;\right\}$$ and, for an index $\mathcal{A}$ and a family of measurable functions $\left\{f_\alpha:\mathbb{R}^d\times\Omega\rightarrow\mathbb{R}^{n_\alpha}\right\}_{\alpha\in\mathcal{A}}$, the sigma algebra generated by the random variables $f_\alpha(x,\omega)$, for $x\in A$ and $\alpha\in\mathcal{A}$, is denoted $$\sigma(f_\alpha(x,\omega)\;|\;x\in A, \alpha\in\mathcal{A}).$$  For domains $U\subset\mathbb{R}^d$, $\USC(U;\mathbb{R}^d)$, $\LSC(U;\mathbb{R}^d)$, $\BUC(U;\mathbb{R}^d)$, $\C(U;\mathbb{R}^d)$, $\Lip(U;\mathbb{R}^d)$, $\C^{0,\beta}(U;\mathbb{R}^d)$ and $\C^k(U;\mathbb{R}^d)$ are the spaces of upper-semicontinuous, lower-semicontinuous, bounded continuous, continuous, Lipschitz continuous, $\beta$-H\"{o}lder continuous and $k$-continuously differentiable functions on $U$ with values in $\mathbb{R}^d$.  Furthermore, $C^\infty_c(\mathbb{R}^d)$ denotes the space of smooth, compactly supported functions on $\mathbb{R}^d$.  The closure and boundary of $U\subset\mathbb{R}^d$ are written $\overline{U}$ and $\partial U$.  For $f:\mathbb{R}^d\rightarrow\mathbb{R}$, the support of $f$ is denoted $\Supp(f)$.  Furthermore, $B_R$ and $B_R(x)$ are respectively the open balls of radius $R$ centered at zero and $x\in\mathbb{R}^d$.  For a real number $r\in\mathbb{R}$, the notation $\left[r\right]$ denotes the largest integer less than or equal to $r$.  Finally, throughout the paper $C$ represents a constant which may change from line to line but is independent of $\omega\in\Omega$ unless otherwise indicated.

\subsection{The Random Environment}

The random environment is indexed by a probability space $(\Omega,\mathcal{F},\mathbb{P})$.  Every element $\omega\in\Omega$ corresponds to an individual realization of the environment as described by the coefficients $A(\cdot,\omega)$ and $b(\cdot,\omega)$ on $\mathbb{R}^d$.  The stationarity of the coefficients is quantified by an \begin{equation}\label{transgroup} \textrm{ergodic group of measure-preserving transformations}\; \left\{\tau_x:\Omega\rightarrow\Omega\right\}_{x\in\mathbb{R}^d}\end{equation} such that the coefficients $A:\mathbb{R}^d\times\Omega\rightarrow\mathcal{S}(d)$ and $b:\mathbb{R}^d\times\Omega\rightarrow\mathbb{R}^d$ are bi-measurable stationary functions satisfying, for each $x,y\in\mathbb{R}^d$ and $\omega\in\Omega$, \begin{equation}\label{stationary} A(x+y,\omega)=A(x,\tau_y\omega)\;\;\textrm{and}\;\;b(x+y,\omega)=b(x,\tau_y\omega).\end{equation}

The diffusion matrix and drift are bounded and Lipschitz uniformly for $\omega\in\Omega$.  There exists $C>0$ such that, for all $x\in\mathbb{R}^d$ and $\omega\in\Omega$,  \begin{equation}\label{bounded} \abs{b(x,\omega)}\leq C\;\;\;\textrm{and}\;\;\;\abs{A(x,\omega)}\leq C, \end{equation} and, for all $x,y\in\mathbb{R}^d$ and $\omega\in\Omega$, \begin{equation}\label{Lipschitz} \abs{b(x,\omega)-b(y,\omega)}\leq C\abs{x-y}\;\;\;\textrm{and}\;\;\;\abs{A(x,\omega)-A(y,\omega)}\leq C\abs{x-y}.\end{equation}  In addition, the diffusion matrix is uniformly elliptic uniformly in $\Omega$.  There exists $\nu>1$ such that, for all $x\in\mathbb{R}^d$ and $\omega\in\Omega$, \begin{equation}\label{elliptic} \frac{1}{\nu} I\leq A(x,\omega)\leq \nu I.\end{equation}

The coefficients satisfy a finite range dependence.  There exists $R>0$ such that, whenever $A,B\subset\mathbb{R}^d$ satisfy $d(A,B)\geq R$, the sigma algebras \begin{equation}\label{finitedep} \sigma(A(x,\omega), b(x,\omega)\;|\;x\in A)\;\;\;\textrm{and}\;\;\; \sigma(A(x,\omega), b(x,\omega)\;|\;x\in B)\;\;\;\textrm{are independent.}\end{equation}  The diffusion matrix and drift satisfy a restricted isotropy condition.  For every orthogonal transformation $r:\mathbb{R}^d\rightarrow\mathbb{R}^d$ which preserves the coordinate axes, for every $x\in\mathbb{R}^d$, \begin{equation}\label{isotropy} (A(rx,\omega),b(rx,\omega))\;\;\;\textrm{and}\;\;\;(rA(x,\omega)r^t,rb(x,\omega))\;\;\;\textrm{have the same law.}\end{equation}  And, the diffusion matrix and drift are a small perturbation of the Laplacian.  There exists $\eta_0>0$, to later be chosen small, such that, for all $x\in\mathbb{R}^d$ and $\omega\in\Omega$, \begin{equation}\label{perturbation} \abs{b(x,\omega)}\leq\eta_0\;\;\textrm{and}\;\;\abs{A(x,\omega)-I}\leq \eta_0.\end{equation}

The final assumptions concern the domain.  The domain \begin{equation}\label{domain_bounded} U\subset\mathbb{R}^d\;\;\textrm{is open and bounded.}\end{equation}  Furthermore, $U$ satisfies an exterior ball condition.  There exists $r_0>0$ so that, for each $x\in\partial U$ there exists $x^*\in\mathbb{R}^d$ such that \begin{equation}\label{exterior} \overline{B}_{r_0}(x^*)\cap \overline{U}=\left\{x\right\}.\end{equation}

To avoid cumbersome statements in what follows, a steady assumption is introduced.  \begin{equation}\label{steady}\textrm{Assume}\;(\ref{transgroup}), (\ref{stationary}), (\ref{bounded}), (\ref{Lipschitz}), (\ref{elliptic}), (\ref{finitedep}), (\ref{isotropy}), (\ref{perturbation}),  (\ref{domain_bounded})\;\textrm{and}\;(\ref{exterior}).\end{equation}

Observe that (\ref{bounded}), (\ref{Lipschitz}) and (\ref{elliptic}) guarantee the well-posedness of the martingale problem associated to to the generator $$\frac{1}{2}\sum_{i,j=1}^da_{ij}(x,\omega)\frac{\partial^2}{\partial x_i\partial x_j}+\sum_{i=1}^db_i(x,\omega)\frac{\partial}{\partial x_i}$$ for each $x\in\mathbb{R}^d$ and $\omega\in\Omega$, see Strook and Varadhan \cite[Chapter~6,7]{SV}.  The corresponding probability measure and expectation on the space of continuous paths $\C([0,\infty);\mathbb{R}^d)$ will be written $P_{x,\omega}$ and $E_{x,\omega}$ where, almost surely with respect to $P_{x,\omega}$, paths $X_t\in\C([0,\infty);\mathbb{R}^d)$ satisfy the stochastic differential equation \begin{equation}\label{sde}\left\{\begin{array}{l} dX_t=b(X_t,\omega)dt+\sigma(X_t,\omega)dB_t, \\ X_0=x,\end{array}\right.\end{equation} for $A(x,\omega)=\sigma(x,\omega)\sigma(x,\omega)^t$, and for $B_t$ some standard Brownian motion under $P_{x,\omega}$ with respect to the canonical right-continuous filtration on $\C([0,\infty);\mathbb{R}^d)$.

The translation and rotational invariance implied in law by (\ref{stationary}) and (\ref{isotropy}) do not imply any invariance properties, in general, for the quenched measures $P_{x,\omega}$.  However, the annealed measures and expectations, defined by the semi-direct products $\mathbb{P}_x=\mathbb{P}\ltimes P_{x,\omega}$ and $\mathbb{E}_x=\mathbb{E}\ltimes E_{x,\omega}$ on $\Omega\times\C([0,\infty);\mathbb{R}^d)$, do satisfy a translational and rotational invariance in the sense that, for all $x,y\in\mathbb{R}^d$, \begin{equation}\label{annealed} \mathbb{E}_{x+y}(X_t)=\mathbb{E}_y(x+X_t)=x+\mathbb{E}_y(X_t),\end{equation} and, for all orthogonal transformations $r$ preserving the coordinate axis and for every $x\in\mathbb{R}^d$, \begin{equation}\label{annealed1} \mathbb{E}_{x}(rX_t)=\mathbb{E}_{rx}(X_t).\end{equation}  This fact plays an important role in \cite{SZ} to preclude, with probability one, the emergence of ballistic behavior of the rescaled process in the asymptotic limit.

Similarly, for each $n\geq 0$ and $x\in\mathbb{R}^d$, let $W^n_x$ denote the Weiner measure on $\C([0,\infty);\mathbb{R}^d)$ and $E^{W^n_x}$ the expectation corresponding to Brownian motion with variance $\alpha_n$ beginning from $x$.  Almost surely with respect to $W^n_x$, paths $X_t\in\C([0,\infty);\mathbb{R}^d)$ satisfy the stochastic differential equation \begin{equation}\label{sde_brownian}\left\{\begin{array}{l} dX_t=\sqrt{\alpha_n}dB_t, \\ X_0=x,\end{array}\right.\end{equation} for $B_t$ some standard Brownian motion under $W^n_x$ with respect to the canonical right-continuous filtration on $\C([0,\infty);\mathbb{R}^d)$.

\subsection{A Remark on Existence and Uniqueness}

The boundedness (\ref{bounded}), Lipschitz continuity (\ref{Lipschitz}) and ellipticity (\ref{elliptic}) of the coefficients together with the boundedness (\ref{domain_bounded}) and regularity (\ref{exterior}) of the domain guarantee the well-posedness, for every $\omega\in\Omega$, of equations like $$\left\{\begin{array}{ll}\frac{1}{2}\tr(A(x,\omega)D^2w)+b(x,\omega)\cdot Dw=g(x) & \textrm{on}\;\;U, \\ u=f(x) & \textrm{on}\;\;\partial U,\end{array}\right.$$ for $f\in C(\partial U)$ and $g\in\C(\overline{U})$, in the class of bounded continuous functions.  See, for instance, Friedman \cite[Chapter~3]{Fr}.  Furthermore, if $\tau$ denotes the exit time from $U$, then $$u(x)=E_{x,\omega}(f(X_\tau)-\int_0^\tau g(X_s)\;dx)\;\;\textrm{on}\;\;\overline{U},$$ see {\O}ksendal \cite[Exercise~9.12]{Oksendal}.

The same assumptions on the coefficients ensure the well-posedness of parabolic equations like $$\left\{\begin{array}{ll} w_t=\frac{1}{2}\tr(A(x,\omega)D^2w)+b(x,\omega)\cdot Dw & \textrm{on}\;\;\mathbb{R}^d\times(0,\infty), \\ w=f(x) & \textrm{on}\;\;\mathbb{R}^d\times\left\{0\right\},\end{array}\right.$$ for continuous initial data $f(x)$ satisfying, for instance and to the extent that it will be applied in this paper, $\abs{f(x)}\leq C(1+\abs{x}^2)$ on $\mathbb{R}^d$, in the class of continuous functions satisfying, locally in time, a quadratic estimate of the same form.  See \cite[Chapter~1]{Fr}.  Furthermore, the solution admits the representation $$w(x,t)=E_{x,\omega}(f(X_t))\;\;\textrm{on}\;\;\mathbb{R}^d\times(0,\infty),$$ see \cite[Exercise~9.12]{Oksendal}.

The analogous formulas hold for the constant coefficient elliptic and parabolic equations associated, for each $n\geq 0$, to the measures $W^n_x$.  Since these facts are well-known, and since the solution to every equation encountered in this paper admits an explicit probabilistic description, the presentation will not further reiterate these points.

\section{The Inductive Framework and Probabilistic Statement}\label{inductive}

In this section, the aspects of \cite{SZ} most relevant to this work will be introduced.  The interested reader will find a full description of the inductive framework in \cite{SZ}, which was later reviewed in the introductions of \cite{F2,F1}.  Forgive, therefore, the terse explanation offered here.

Fix the dimension \begin{equation}\label{dimension} d\geq 3, \end{equation} and fix a H\"older exponent \begin{equation}\label{Holderexponent} \beta\in\left(0,\frac{1}{2}\right]\;\;\textrm{and a constant}\;\;a\in \left(0,\frac{\beta}{1000d}\right]. \end{equation}

Let $L_0$ be a large integer multiple of five.  For each $n\geq 0$, inductively define \begin{equation}\label{L} \ell_n=5\left[\frac{L_n^a}{5}\right]\;\;\textrm{and}\;\;L_{n+1}=\ell_n L_n, \end{equation} so that, for $L_0$ sufficiently large, it follows that $\frac{1}{2}L_n^{1+a}\leq L_{n+1}\leq 2L_n^{1+a}$.  For each $n\geq 0$, for $c_0>0$, let \begin{equation}\label{kappa} \kappa_n=\exp(c_0(\log\log(L_n))^2)\;\;\textrm{and}\;\;\tilde{\kappa}_n=\exp(2c_0(\log\log(L_n))^2),\end{equation} where, as $n$ tends to infinity, notice that $\kappa_n$ is eventually dominated by every positive power of $L_n$.  Furthermore, define, for each $n\geq 0$, \begin{equation}\label{D} D_n=L_n\kappa_n\;\;\textrm{and}\;\;\tilde{D}_n=L_n\tilde{\kappa}_n,\end{equation} where the preceding remark indicates the scales $D_n$ and $\tilde{D}_n$ are larger but grow comparably with the previously defined scales $L_n$.

The following constants enter into the probabilistic statements below.  Fix $m_0\geq 2$ satisfying \begin{equation}\label{m0} (1+a)^{m_0-2}\leq 100<(1+a)^{m_0-1}, \end{equation}  and $\delta>0$ and $M_0>0$ satisfying \begin{equation}\label{delta} \delta=\frac{5}{32}\beta\;\;\textrm{and}\;\;M_0\geq100d(1+a)^{m_0+2}.\end{equation}  In the arguments to follow, it will be essential that these assumptions guarantee $\delta$ and $M_0$ are sufficiently larger than $a$.

In order to apply the finite range dependence, it will be frequently necessary to introduce a stopped version of the process.  Define for every element $X_t\in \C([0,\infty);\mathbb{R}^d)$ the path \begin{equation}\label{prob_tail}X_t^*=\sup_{0\leq s\leq t}\abs{X_s-X_0},\end{equation} and, for each $n\geq 0$, the stopping time $$T_n=\inf\left\{s\geq 0\;|\;X_s^*\geq \tilde{D}_n\right\}.$$  The effective diffusivity of the ensemble at scale $L_n$ is defined by $$\alpha_n=\frac{1}{2d}\mathbb{E}_0\left(\abs{X_{L_n^2\wedge T_n}}^2\right),$$ where the localization is applied in order to exploit the diffusion's mixing properties.  The convergence of the $\alpha_n$ to a limiting diffusivity $\overline{\alpha}$ is proven in \cite[Proposition~5.7]{SZ}.

\begin{thm}\label{effectivediffusivity} Assume (\ref{steady}).  There exists $L_0$ and $c_0$ sufficiently large and $\eta_0>0$ sufficiently small such that, for all $n\geq 0$, $$\frac{1}{2\nu}\leq \alpha_n\leq 2\nu\;\;\textrm{and}\;\;\abs{\alpha_{n+1}-\alpha_n}\leq L_n^{-(1+\frac{9}{10})\delta},$$  which implies the existence of $\overline{\alpha}>0$ satisfying $$\frac{1}{2\nu}\leq \overline{\alpha}\leq 2\nu\;\;\textrm{and}\;\;\lim_{n\rightarrow\infty}\alpha_n=\overline{\alpha}.$$\end{thm}

The results of \cite{SZ} obtain an effective comparison on the parabolic scale $(L_n, L_n^2)$ in space and time, with improving probability as $n\rightarrow\infty$, between the solutions \begin{equation}\label{prob_eq} \left\{\begin{array}{ll} u_t=\frac{1}{2}\tr(A(x,\omega)D^2u)+b(x,\omega)\cdot Du & \textrm{on}\;\;\mathbb{R}^d\times(0,\infty), \\ u=f(x) & \textrm{on}\;\;\mathbb{R}^d\times\left\{0\right\},\end{array}\right.\end{equation} and solutions to the approximate limiting equation \begin{equation}\label{prob_approx}\left\{\begin{array}{ll} u_{n,t}=\frac{\alpha_n}{2}\Delta u_n & \textrm{on}\;\;\mathbb{R}^d\times(0,\infty), \\ u_n=f(x) & \textrm{on}\;\;\mathbb{R}^d\times\left\{0\right\}.\end{array}\right.\end{equation}  In order to simplify the notation define, for each $n\geq 0$, the operators \begin{equation}\label{prob_operators} R_nf(x)=u(x,L_n^2)\;\;\textrm{and}\;\;\overline{R}_nf(x)=u_n(x,L_n^2),\end{equation} and the difference operator \begin{equation}\label{prob_difference} S_nf(x)=R_nf(x)-\overline{R}_nf(x).\end{equation}

Since solutions of (\ref{prob_eq}) will not, in general, be effectively comparable with solutions of (\ref{prob_approx}) globally in space, it is necessary to introduce a cutoff function.   For each $v>0$, let \begin{equation}\label{cutoff} \chi(y)=1\wedge(2-\abs{y})_+\;\;\textrm{and}\;\;\chi_{v}(y)=\chi\left(\frac{y}{v}\right), \end{equation}  and define, for each $x\in\mathbb{R}^d$ and $n\geq 0$, \begin{equation}\label{cutoff1}  \chi_{n,x}(y)=\chi_{30\sqrt{d}L_n}(y-x).\end{equation}  Furthermore, in order to account for the scaling of the initial data which appears in (\ref{intro_eq}), the comparison of the solutions is necessarily obtained with respect to the rescaled global H\"older-norms, defined for each $n\geq 0$, \begin{equation}\label{prob_Holder} \abs{f}_n=\norm{f}_{L^\infty(\mathbb{R}^d)}+\sup_{x\neq y}L_n^\beta\frac{\abs{f(x)-f(y)}}{\abs{x-y}^\beta}.\end{equation}  See for instance the introductions of \cite{F1,SZ} for a more complete discussion concerning the necessity of these norms as opposed, say, to attempting an (in general, false) $L^\infty$-comparison.

The following control is the statement propagated by the arguments of \cite{SZ}, and expresses the desired comparison between solutions (\ref{prob_eq}) and (\ref{prob_approx}), as written using the operator $S_n$ from (\ref{prob_difference}) and localized by $\chi_{n,x}$ from (\ref{cutoff1}), in terms of the $\abs{\cdot}_n$-norm from (\ref{prob_Holder}) of the initial data.

Note carefully that this statement is not true, in general, for all triples $x\in\mathbb{R}^d$, $\omega\in\Omega$ and $n\geq 0$.  However, as described below, it is shown in \cite[Proposition~5.1]{SZ} that such controls are available for large $n$, with high probability, on a large portion of space.

\begin{con}\label{Holder}  Fix $x\in\mathbb{R}^d$, $\omega\in\Omega$ and $n\geq 0$.  Then, for each $f\in C^{0,\beta}(\mathbb{R}^d)$, $$\abs{\chi_{n,x}S_nf}_n\leq L_n^{-\delta}\abs{f}_n.$$\end{con}

It will also be necessary to obtain tail-estimates for the diffusion in random environment.  Recalling that $P_{x,\omega}$ is the measure on the space of continuous paths describing the diffusion beginning from $x\in\mathbb{R}^d$ and associated to the generator $$\frac{1}{2}\sum_{i,j=1}^da_{ij}(x,\omega)\frac{\partial^2}{\partial x_i \partial x_j}+\sum_{i=1}^db_i(x,\omega)\frac{\partial}{\partial x_i},$$ the type of control propagated in \cite{SZ} involves exponential estimates for the probability under $P_{x,\omega}$ that the maximal excursion $X^*_{L_n^2}$ defined in (\ref{prob_tail}) is large with respect to the time elapsed.

As with Control \ref{Holder}, it is simply not true in general that this type of estimate is satisfied for all triples $(x,\omega,n)$.  However, it is shown in \cite[Proposition~2.2]{SZ} that such controls are available for large $n$, with high probability, on a large portion of space.

\begin{con}\label{localization}  Fix $x\in\mathbb{R}^d$, $\omega\in\Omega$ and $n\geq 0$.  For each $v\geq D_n$, for all $\abs{y-x}\leq 30\sqrt{d}L_n$, $$P_{y,\omega}(X^*_{L_n^2}\geq v)\leq \exp(-\frac{v}{D_n}).$$\end{con}

It is necessary to obtain a lower bound in probability for the event, defined for each $n\geq 0$ and $x\in\mathbb{R}^d$, \begin{equation}\label{mainevent} B_n(x)=\left\{\;\omega\in\Omega\;|\;\textrm{Controls \ref{Holder} and \ref{localization} hold for the triple}\;(x,\omega,n).\;\right\}.\end{equation}  Notice that, in view of (\ref{stationary}), for all $x\in\mathbb{R}^d$ and $n\geq 0$, \begin{equation}\label{mainevent1}\mathbb{P}(B_n(x))=\mathbb{P}(B_n(0)),\end{equation} and observe that $B_n(0)$ does not include the control of traps described in \cite[Proposition~3.3]{SZ}, which play in important role in propagating Control \ref{Holder}, and from which the arguments of this paper have no further need.

The following theorem proves that the compliment of $B_n(0)$ approaches zero as $n$ tends to infinity, see \cite[Theorem~1.1]{SZ}.

\begin{thm}\label{induction}  Assume (\ref{steady}).  There exist $L_0$ and $c_0$ sufficiently large and $\eta_0>0$ sufficiently small such that, for each $n\geq 0$, $$\mathbb{P}\left(\Omega\setminus B_n(0)\right)\leq L_n^{-M_0}.$$\end{thm}

Henceforth, the constants $L_0$, $c_0$ and $\eta_0$ are fixed to satisfy the hypothesis of Theorems \ref{effectivediffusivity} and \ref{induction} appearing above.  \begin{equation}\label{constants} \textrm{Fix constants}\;L_0, c_0\;\textrm{and}\;\eta_0\;\textrm{satisfying the hypothesis of Theorems \ref{effectivediffusivity} and \ref{induction}.}\end{equation}

The events which, following an application of the Borel-Cantelli lemma, come to define the event on which Theorem \ref{intro_main} is obtained are chosen to ensure that Controls \ref{Holder} and \ref{localization} are satisfied for a sufficiently small scale as compared with $\frac{1}{\epsilon}$.  Fix the smallest integer $\overline{m}>0$ satisfying the inequality \begin{equation}\label{prob_m} \overline{m}>1-\frac{\log(1-2a-a^2)}{\log(1+a)},\end{equation} where the definition of $L_n$ in (\ref{L}) implies that $L_{n+1}L_{n-\overline{m}}\leq L_{n-1}^2$ for all $n\geq0$ sufficiently large.

The idea will be to use Theorem \ref{induction} in order to obtain Controls \ref{Holder} and \ref{localization} at scale $L_{n-\overline{m}}$ on the entirety of the rescaled domain $U/\epsilon$ whenever $L_n\leq \frac{1}{\epsilon}<L_{n+1}$.  Since, for all $n\geq 0$ sufficiently large, it follows from the boundedness of $U$ and (\ref{L}) that, whenever $L_n\leq\frac{1}{\epsilon}<L_{n+1}$, the rescaled domain $U/\epsilon$ is contained in what becomes the considerably larger set $[-\frac{1}{2}L_{n+2}^2, \frac{1}{2}L_{n+2}^2]^d$, define, for each $n\geq \overline{m}$, \begin{multline}\label{prob_event_1} A_n=\left\{\;\omega\in\Omega\;|\;\omega\in B_m(x)\;\;\textrm{for all}\;\;x\in L_m\mathbb{Z}^d\cap[-L_{n+2}^2, L_{n+2}^2]^d\;\;\textrm{and}\right. \\ \left.\textrm{for all}\;\;n-\overline{m}\leq m\leq n+2\right\}.\end{multline}  The following proposition proves that, as $n\rightarrow\infty$, the probability of the events $A_n$ rapidly approaches one, since the exponent $$2d(1+a)^2-\frac{M_0}{2}<0$$ is negative owing to (\ref{Holderexponent}) and (\ref{delta}).

\begin{prop}\label{prob_probability}  Assume (\ref{steady}) and (\ref{constants}).  For each $n\geq \overline{m}$, for $C>0$ independent of $n$, $$\mathbb{P}(\Omega\setminus A_n)\leq CL_n^{2d(1+a)^2-\frac{1}{2}M_0}.$$\end{prop}

\begin{proof}  Fix $n\geq\overline{m}$.  Theorem \ref{induction} and (\ref{mainevent1}) imply using (\ref{L}) that, for $C>0$ independent of $n$, $$\mathbb{P}(\Omega\setminus A_n)\leq \sum_{m=n-\overline{m}}^{n+2}(\frac{L_{n+2}^2}{L_{m}})^dL_m^{-M_0}\leq C\sum_{m=n-\overline{m}}^{n+2}L_n^{2d(1+a)^2-2d(1+a)^{m-n}-M_0(1+a)^{m-n}}.$$  Therefore, $$\mathbb{P}(\Omega\setminus A_n)\leq CL_n^{2d(1+a)^2-M_0(1+a)^{-\overline{m}}},$$ which, since the definition of $\overline{m}$ implies that $$(1+a)^{-\overline{m}}\geq (1+a)(\frac{2}{1-a}-(1+a))\geq\frac{1}{2},\;\;\textrm{yields}\;\;\mathbb{P}(\Omega\setminus A_n)\leq CL_n^{2d(1+a)^2-\frac{M_0}{2}}$$ and completes the proof.\end{proof}

\section{An Upper Bound for the Exit Time of the Process in Random Environment}\label{section_exit_time}

The purpose of this section is to obtain an upper bound in probability for the exit time from the rescaled domain $U/\epsilon$ of the process associated to the generator \begin{equation}\label{exit_generator}\frac{1}{2}\sum_{i,j=1}^da_{ij}(x,\omega)\frac{\partial^2}{\partial x_i \partial x_j}+\sum_{k=1}^db_i(x,\omega)\frac{\partial}{\partial x_i}.\end{equation}  The reason for obtaining such an estimate will be seen in Section 5, where the process in random environment is coupled with high probability to a deterministic Brownian motion.  Since this coupling cannot be expected to hold globally in time, it is necessary to ensure with high probability that the exit time from $U/\epsilon$ occurs before the estimates deteriorate.

It will be shown that, as a consequence of the H\"older estimate stated in Control \ref{Holder}, whenever the environment and scale satisfy $\omega\in A_n$ and $L_n\leq \frac{1}{\epsilon}<L_{n+1}$ then, as $n\rightarrow\infty$, the exit time from the rescaled domain $U/\epsilon$ occurs before time $L_{n+2}^2$ with overwhelming probability.  Define, for each $\epsilon>0$, the $\C([0,\infty);\mathbb{R}^d)$ exit time \begin{equation}\label{exit_time}\tau^\epsilon=\inf\left\{\;t\geq 0\;|\;X_t\notin U/\epsilon\;\right\}=\inf\left\{\;t\geq 0\;|\;\epsilon X_t\notin U\;\right\},\end{equation}  where the final equality is particularly prescient in view of (\ref{intro_sz}) and the scaling associated to the generator $$\frac{1}{2}\sum_{i,j=1}^da_{ij}(\frac{x}{\epsilon},\omega)\frac{\partial^2}{\partial x_i \partial x_j}+\frac{1}{\epsilon}\sum_{i=1}^db_i(\frac{x}{\epsilon},\omega)\frac{\partial}{\partial x_i}.$$  In terms of this rescaled generator, the following proposition proves that, for environments $\omega\in A_n$ and scales $L_n\leq\frac{1}{\epsilon}<L_{n+1}$, as $n\rightarrow\infty$, paths $\epsilon X_{t/\epsilon^2}$ exit $U$ with overwhelming probability prior to time $\epsilon^2L_{n+2}^2$.

\begin{prop}\label{exit_main}  Assume (\ref{steady}) and (\ref{constants}).  For all $n$ sufficiently large, for every $\omega\in A_n$, for all $\epsilon>0$ satisfying $L_n\leq \frac{1}{\epsilon}<L_{n+1}$, for $C>0$ independent of $n$, $$\sup_{x\in \overline{U}}P_{\frac{x}{\epsilon},\omega}(\tau^\epsilon>L_{n+2}^2)\leq C L_n^{-da}.$$\end{prop}

\begin{proof}  Using the boundedness of the domain in (\ref{domain_bounded}), choose $R\geq 1$ satisfying $\overline{U}\subset B_R$ and choose $n_1\geq 0$ such that, for every $n\geq n_1$, \begin{equation}\label{exit_main_1}L_{n+1}\overline{U}\subset L_{n+1}B_R\subset[-L_{n+2}^2, L_{n+2}^2]^d.\end{equation} Henceforth, fix $n\geq n_1$, $\omega\in A_n$ and $L_n\leq \frac{1}{\epsilon}<L_{n+1}$.

Define the smooth cutoff function satisfying $0\leq \chi_{B_R}\leq 1$ with $$\chi_{B_{R}}(x)=\left\{\begin{array}{ll} 1 & \textrm{if}\;\;x\in \overline{B}_R, \\ 0 & \textrm{if}\;\;x\in \mathbb{R}^d\setminus B_{R+1},\end{array}\right.$$ and observe that, for a constant $C>0$ independent of $\epsilon>0$, since $L_n\leq \frac{1}{\epsilon}<L_{n+1}$, \begin{equation}\label{exit_main_2} \abs{\chi_{B_R}(\epsilon x)}_{n+2}\leq 1+C\frac{L_{n+2}}{L_n}\leq CL_n^{2a+a^2}.\end{equation}  Then, consider solutions $$\left\{\begin{array}{ll} v^\epsilon_t=\frac{1}{2}\tr(A(x,\omega)D^2v^\epsilon)+b(x,\omega)\cdot Dv^\epsilon & \textrm{on}\;\;\mathbb{R}^d\times(0,\infty), \\ v^\epsilon=\chi_{B_R}(\epsilon x) & \textrm{on}\;\;\mathbb{R}^d\times\left\{0\right\},\end{array}\right.$$ which admit the representation $$v^\epsilon(x,t)=E_{x,\omega}(\chi_{B_R}(\epsilon X_t))\geq P_{x,\omega}(\epsilon X_t\in B_R)\geq P_{x,\omega}(\epsilon X_t \in U).$$ Therefore, \begin{equation}\label{exit_main_4}1-v^\epsilon(x,t)\leq P_{x,\omega}(\epsilon X_t\notin U)\leq P_{x,\omega}(\tau^\epsilon\leq t).\end{equation}

The function $v^\epsilon$ will be compared via Control \ref{Holder} with the solution $$\left\{\begin{array}{ll} \overline{v}^\epsilon_t=\frac{\alpha_{n+2}}{2}\Delta \overline{v}^\epsilon & \textrm{on}\;\;\mathbb{R}^d\times(0,\infty), \\ \overline{v}^\epsilon=\chi_{B_R}(\epsilon x) & \textrm{on}\;\;\mathbb{R}^d\times\left\{0\right\}.\end{array}\right.$$  The conditions $\omega\in A_n$ and (\ref{exit_main_1}) guarantee that for every $x\in \overline{U}/\epsilon$ the conclusion of Control \ref{Holder} is satisfied and, therefore, using (\ref{Holderexponent}), (\ref{L}) and (\ref{exit_main_2}), for $C>0$ independent of $n$, \begin{equation}\label{exit_main_3} \sup_{x\in \overline{U}}\abs{v^\epsilon(x,L_{n+2}^2)-\overline{v}^\epsilon(x,L_{n+2}^2)}\leq CL_{n+2}^{-\delta}L_n^{2a+a^2}\leq CL_n^{2a+a^2-\delta(1+a)^2}\leq CL_n^{3a-\delta}.\end{equation}

To conclude, the size of $\overline{v}^\epsilon(x,L_{n+2}^2)$, which measures the likelihood that a Brownian motion with variance $\alpha_{n+2}$ and beginning from $x$ resides outside $B_{\frac{R+1}{\epsilon}}$ at time $L_{n+2}^2$, is estimated using Theorem \ref{effectivediffusivity} and the Green's function.  For each $x\in \overline{U}/\epsilon$, since $\frac{1}{\epsilon}< L_{n+1}$ and $R\geq 1$, for $C>0$ independent of $n$, $$\overline{v}^\epsilon(x,L_{n+2}^2)\leq \int_{B_{\frac{4R}{\epsilon}}(x)}(4\pi \alpha_{n+2}L_{n+2}^2)^{-\frac{d}{2}}\exp(-\frac{\abs{y-x}^2}{4\alpha_{n+2} L_{n+2}^2})\;dy\leq C(\epsilon L_{n+2})^{-d}\leq CL_n^{-da(1+a)}.$$  Therefore, in view of (\ref{exit_main_3}), for each $x\in \overline{U}/\epsilon$, for $C>0$ independent of $n$, \begin{equation}\label{exit_main_5} 1-v^\epsilon(x,L_{n+2}^2)\geq 1-\overline{v}^\epsilon(x,L_{n+2}^2)-\abs{v^\epsilon(x,L_{n+2}^2)-\overline{v}^\epsilon(x,L_{n+2}^2)}\geq 1-CL_n^{-da(1+a)}-CL_n^{3a-\delta}.\end{equation}  So, using (\ref{exit_main_4}), since (\ref{Holderexponent}) and (\ref{delta}) imply $da<da(1+a)<\delta-3a$, for $C>0$ independent of $n$, $$\sup_{x\in \overline{U}}P_{\frac{x}{\epsilon},\omega}(\tau^\epsilon>L_{n+2}^2)\leq CL_n^{-da},$$ which completes the argument.  \end{proof}

\section{The Global Coupling}\label{coupling}

The comparison implied by Control \ref{Holder} on scale $(L_n, L_n^2)$ between the vector-valued solutions of the parabolic equation \begin{equation}\label{couple_eq}\left\{\begin{array}{ll} u_t=\frac{1}{2}\tr(A(x,\omega)D^2u)+b(x,\omega)\cdot Du & \textrm{on}\;\;\mathbb{R}^d\times(0,\infty), \\ u=\frac{x}{L_n} & \textrm{on}\;\;\mathbb{R}^d\times\left\{0\right\},\end{array}\right.\end{equation} and the approximate homogenized equation \begin{equation}\label{couple_approx} \left\{\begin{array}{ll} u_{n,t}=\frac{\alpha_n}{2}\Delta u_n & \textrm{on}\;\;\mathbb{R}^d\times(0,\infty), \\ u_n=\frac{x}{L_n} & \textrm{on}\;\;\mathbb{R}^d\times\left\{0\right\},\end{array}\right.\end{equation} asserts that, after using the localization estimate implied by Control \ref{localization} and the choice of constants (\ref{L}), (\ref{kappa}) and (\ref{D}) to localize and bound the initial data with respect to the $\abs{\cdot}_n$-norm, \begin{equation}\label{couple_approx_100}\abs{u(0,L_n^2)-u_n(0,L_n^2)}=\abs{E_{0,\omega}(\frac{1}{L_n}X_{L_n^2})-E^{W_0^n}(\frac{1}{L_n}X_{L_n^2}))}\leq C\tilde{\kappa}_n L_n^{-\delta},\end{equation} where $W^n_x$ the Weiner measure on $\C([0,\infty);\mathbb{R}^d)$ corresponding to Brownian motion with variance $\alpha_n$ beginning from $x$.

Very formally, then, provided (what will be discrete) copies of the diffusion in random environment $\tilde{X}_t$ and Brownian motion $\tilde{B}_t$ are chosen carefully and are defined with respect to the same measure on an auxiliary probability space $(\tilde{\Omega},\tilde{\mathcal{F}},\tilde{\mathbb{P}})$, a Chebyshev inequality should yield $$(\frac{\gamma}{L_n})^\beta\tilde{\mathbb{P}}(\abs{\tilde{X}_{L_n^2}-\tilde{B}_{L_n^2}}^\beta\geq \gamma^\beta)\leq CL_n^{-\delta}\tilde{\kappa}_n,$$ which implies \begin{equation}\label{couple_goal} \tilde{\mathbb{P}}(\abs{\tilde{X}_{L_n^2}-\tilde{B}_{L_n^2}}\geq \gamma)\leq CL_n^{-\delta}\tilde{\kappa}_n(\frac{L_n}{\gamma})^\beta.\end{equation}  The purpose of this section will be to formalize and iterate this intuition along a discrete sequence of time steps, where the work comes in constructing the processes $\tilde{X}_t$ and $\tilde{B}_t$ such that the integration and absolute value in the version of (\ref{couple_approx_100}) with respect to $\tilde{\mathbb{P}}$ is commuted, and thereby justifies truly the application of the Chebyshev inequality.

Solutions of (\ref{couple_eq}) with initial data $f(x)$ admit a representation in terms of the Green's function $$p_{t,\omega}(x,y):[0,\infty)\times\mathbb{R}^d\times\mathbb{R}^d\rightarrow\mathbb{R},$$ which is the density of the diffusion beginning from $x$ in environment $\omega$ at time $t$.  See \cite[Chapter~1]{Fr} for a detailed discussion of the existence and regularity of these densities, and which follow from assumptions (\ref{bounded}), (\ref{Lipschitz}) and (\ref{elliptic}).  The formula for the solution is then $$u(x,t)=E_{x,\omega}(f(X_t))=\int_{\mathbb{R}^d}p_{t,\omega}(x,y)f(y)\;dy.$$  Similarly, solutions of (\ref{couple_approx}) with initial data $f(x)$ admit the analogous representation in terms of the heat kernel $$\overline{u}_n(x,t)=E^{W_x^n}(f(X_t))=\int_{\mathbb{R}^d}(4\pi\alpha_n t)^{-\frac{d}{2}}\exp(-\frac{\abs{y-x}^2}{4\alpha_n t})f(y)\;dy.$$

To simplify the notation in what follows, for each $n\geq 0$, define $$p_{n,\omega}(x,y)=p_{L_n^2,\omega}(x,y),$$ and the analogous heat kernel $$\overline{p}_n(x,y)=(4\pi\alpha_n L_n^2)^{-\frac{d}{2}}\exp(-\frac{\abs{y-x}^2}{4\alpha_n L_n^2}).$$  The following proposition constructs a Markov process $(X_k,\overline{X}_k)$ on the space $(\mathbb{R}^d\times\mathbb{R}^d)^{\mathbb{N}}$ such that the transition probabilities of first coordinate $X_k$ are determined by $p_{n,\omega}$ and, those of the second coordinate $\overline{X}_k$ by $\overline{p}_n$.  Furthermore, the difference $\abs{X_k-\overline{X_k}}$ satisfies a version of (\ref{couple_goal}) with respect to the underlying measure.

The construction follows closely the proof of \cite[Proposition~3.1]{SZ}, and is included for the convenience of the reader and due to the mildly different formulation adapted to the arguments in this paper.  The proof relies upon the Kantorovich-Rubinstein Theorem, see Dudley \cite[Theorem~11.8.2]{D}, applied to the metrics $$d_n(x,y)=\abs{\frac{x-y}{L_n}}^\beta,$$ where $0<\beta\leq\frac{1}{2}$ was fixed in (\ref{Holderexponent}).  The theorem states that every pair of probability measures $\nu$ and $\nu'$ on $\mathbb{R}^d$ assigning finite mass to the metric $d_n$, in the sense that \begin{equation}\label{couple_kr_1} \int_{\mathbb{R}^d}d_n(x,0)\;\nu(dx)<\infty\;\;\textrm{and}\;\;\int_{\mathbb{R}^d}d_n(x,0)\;\nu'(dx)<\infty,\end{equation} satisfy the equality \begin{multline}\label{couple_kr} D_n(\nu,\nu')=\sup\left\{\abs{\int f\;d\nu-\int f\;d\nu'}\;|\;\abs{f(x)-f(y)}\leq d_n(x,y)\;\;\textrm{on}\;\;\mathbb{R}^d\times\mathbb{R}^d\right\} \\ =\inf\left\{\int_{\mathbb{R}^d\times\mathbb{R}^d}d_n(x,x')\;\rho(dx,dx')\;|\;\rho\;\textrm{is a probability measure on}\;\mathbb{R}^d\times\mathbb{R}^d\right. \\ \left. \textrm{with first and second marginals}\;\nu\;\textrm{and}\;\nu'\right\}.\end{multline}  The function $D_n(\cdot,\cdot)$ is sometimes referred to as the Kantorovich-Rubinstein or Wasserstein metric.

The choice of constants in the following proposition will be applied to spacial scales satisfying $L_{n}\leq \frac{1}{\epsilon}<L_{n+1}$.  Therefore, in view of Proposition \ref{exit_main}, the coupling remains effective up to and past a point that the diffusion has exited the domain with overwhelming probability.  Furthermore, the time steps $L_{n-\overline{m}}^2$ are chosen to be much smaller than the scale $\frac{1}{\epsilon^2}$ in time as determined by the definition of $\overline{m}$ in (\ref{prob_m}).

\begin{prop}\label{couple_main}  Assume (\ref{steady}) and (\ref{constants}).  For every $\omega\in\Omega$, for every $x\in\mathbb{R}^d$, there exists a measure $Q_{n,x}$ on the canonical sigma algebra of the space $(\mathbb{R}^d\times\mathbb{R}^d)^{\mathbb{N}}$ such that, under $Q_{n,x}$, the coordinate processes $X_k$ and $\overline{X}_k$ respectively have the law of a Markov chain on $\mathbb{R}^d$, starting from $x$, with transition kernels $p_{n-\overline{m},\omega}(\cdot,\cdot)$ and $\overline{p}_{n-\overline{m}}(\cdot,\cdot)$.

Furthermore, for every $n\geq \overline{m}$, $\omega\in A_n$ and $x\in[-\frac{1}{2}L_{n+2}^2, \frac{1}{2}L_{n+2}^2]^d$, for $C>0$ independent of $n$, \begin{equation}\label{couple_main_0}Q_{n,x}(\abs{X_k-\overline{X}_k}\geq \gamma\;|\;\textrm{for some}\;0\leq k\leq 2(\frac{L_{n+2}}{L_{n-\overline{m}}})^2)\leq C(\frac{L_{n-\overline{m}}}{\gamma})^\beta(\frac{L_{n+2}}{L_{n-\overline{m}}})^4\tilde{\kappa}_{n-\overline{m}}L_{n-\overline{m}}^{-\delta}.\end{equation}\end{prop}

\begin{proof}  Fix $n\geq \overline{m}$ and $\omega\in \Omega$.  Let $M_1(\mathbb{R}^d\times\mathbb{R}^d)$ denote the set of probability measures on $\mathbb{R}^d\times\mathbb{R}^d$ with the topology of weak convergence.  Exponential estimates imply, for each $x\in\mathbb{R}^d$, the integrals in (\ref{couple_kr_1}) corresponding to the kernels $\nu_x=p_{n-\overline{m}}(x,\cdot)$ and $\nu'_x=\overline{p}_{n-\overline{m}}(x,\cdot)$ are finite, see \cite[Chapter~1,Theorem~12]{Fr}.  The Kantorovich-Rubinstein theorem, see (\ref{couple_kr}), therefore implies that, for each $x\in\mathbb{R}^d$, the subset \begin{multline}\label{couple_main_1}K_x=\left\{\;\rho\in M_1(\mathbb{R}^d\times\mathbb{R}^d)\;|\;\rho\;\;\textrm{has marginals}\;p_{n-\overline{m},\omega}(x,\cdot)\;\textrm{and}\;\overline{p}_{n-\overline{m}}(x,\cdot),\right. \\ \left. \textrm{and}\;D_{n-\overline{m}}(p_{n-\overline{m},\omega}(x,\cdot),\overline{p}_{n-\overline{m}}(x,\cdot))=\int_{\mathbb{R}^d\times\mathbb{R}^d}d_n(x_1,x_2)\;\rho(dx_1,dx_2)\;\;\right\}\end{multline}  is non-empty and compact.

Furthermore, if $\left\{x_i\in\mathbb{R}^d\right\}_{i=1}^\infty$ is a sequence with limit $x_\infty\in\mathbb{R}^d$ then any sequence $\left\{\rho_i\in K_{x_i}\right\}_{i=1}^\infty$ is tight and has limit $\rho_\infty$ satisfying \begin{multline}\label{gc_main_2}\int_{\mathbb{R}^d\times\mathbb{R}^d}d_{n-\overline{m}}(x_1,x_2)\;\rho_{\infty}(dx_1,dx_2)\leq \\ \liminf_{i\rightarrow\infty}\int_{\mathbb{R}^d\times\mathbb{R}^d}d_{n-\overline{m}}(x_1,x_2)\;\rho_i(dx_1,dx_2)=D_{n-\overline{m}}(p_{n-\overline{m},\omega}(x_\infty,\cdot),\overline{p}_{n-\overline{m}}(x_\infty,\cdot)),\end{multline} where the final inequality follows using (\ref{couple_kr}), (\ref{couple_main_1}) and the triangle inequality satisfied by $D_{n-\overline{m}}$.  Since (\ref{gc_main_2}) implies that $\rho_\infty\in K_{z_\infty}$, using \cite[Lemma~12.1.8, Theorem~12.1.10]{SV}, there exists a measurable map from $\mathbb{R}^d$ to $M_1(\mathbb{R}^d\times\mathbb{R}^d)$ satisfying \begin{equation}\label{gc_main_3}x\rightarrow \tilde{\rho}_x\in M_1(\mathbb{R}^d\times\mathbb{R}^d)\;\;\textrm{with}\;\;\tilde{\rho}_x\in K_x.\end{equation}

The transition kernel of the Markov chain beginning at $(x,y)\in(\mathbb{R}^d\times\mathbb{R}^d)$ is denoted $\tilde{p}_{x,y}\in M_1(\mathbb{R}^d\times\mathbb{R}^d)$ and is defined, for each $g\in L^\infty(\mathbb{R}^d\times\mathbb{R}^d)$, by the relation  \begin{equation}\label{gc_main_4} \int_{\mathbb{R}^d\times\mathbb{R}^d}g(x_1,x_2)\;\tilde{p}_{x,y}(dx_1,dx_2)=\int_{\mathbb{R}^d\times\mathbb{R}^d}g(x_1, x_2-x+y)\;\tilde{\rho}_x(dx_1,dx_2).\end{equation}  For each $x\in\mathbb{R}^d$, the measure $Q_{n,x}$ is defined as the law of the Markov chain $(X_k,\overline{X}_k)$ with transition kernel $\tilde{p}_{\cdot,\cdot}$ and and initial distribution $(x,x)$.

Notice that, if $A\subset\mathbb{R}^d$ is a Borel subset and $k\geq 0$, then, using (\ref{couple_main_1}), (\ref{gc_main_2}) and (\ref{gc_main_4}), for each $x\in\mathbb{R}^d$ and $(y,z)\in\mathbb{R}^d\times\mathbb{R}^d$, $$Q_{n,x}(X_{k+1}\in A\;|\;(X_k,\overline{X}_k)=(y,z))=\int_{A\times\mathbb{R}^d}\tilde{\rho}_y(dx_1,dx_2)=\int_Ap_{n-\overline{m},\omega}(y,x_1)\;dx_1,$$ and, similarly, \begin{multline*}  Q_{n,x}(\overline{X}_{k+1}\in A\;|\;(X_k,\overline{X}_k)=(y,z))=\int_{\mathbb{R}^d\times(A+y-z)}\tilde{\rho}_y(dx_1,dx_2) \\ =\int_{A+y-z}\overline{p}_{n-\overline{m}}(y,x_2)\;dx_2=\int_A\overline{p}_{n-\overline{m}}(z,x_2)\;dx_2,\end{multline*} where the final line uses the translation invariance and symmetry of the heat kernel.  This completes the proof of existence.  It remains to show (\ref{couple_main_0}).

Fix $n\geq \overline{m}$, $\omega\in A_n$ and $x\in[-\frac{1}{2}L_{n+2}^2, \frac{1}{2}L_{n+2}^2]^d$.  Let $0\leq k\leq 2(\frac{L_{n+2}^2}{L_{n-\overline{m}}^2})$ be arbitrary.  The triangle inequality and definition of $d_{n-\overline{m}}$ imply that, writing $E^{Q_{n,x}}$ for the expectation with respect to $Q_{n,x}$, $$E^{Q_{n,x}}(d_{n-\overline{m}}(X_k,\overline{X}_k))\leq E^{Q_{n,x}}(d_{n-\overline{m}}(X_{k-1},\overline{X}_{k-1}))+E^{Q_{n,x}}(d_{n-\overline{m}}(X_k, \overline{X}_k-\overline{X}_{k-1}+X_{k-1})),$$ where, using (\ref{couple_main_1}), (\ref{gc_main_3}), (\ref{gc_main_4}) and the strong Markov property, $$E^{Q_{n,x}}(d_{n-\overline{m}}(X_k, \overline{X}_k-\overline{X}_{k-1}+X_{k-1}))=E^{Q_{n,x}}(D_{n-\overline{m}}(p_{n-\overline{m},\omega}(X_{k-1},\cdot), \overline{p}_{n-\overline{m}}(X_{k-1},\cdot))).$$  Therefore, \begin{multline}\label{couple_main_3} E^{Q_{n,x}}(d_{n-\overline{m}}(X_k,\overline{X}_k))\leq \\ E^{Q_{n,x}}(d_{n-\overline{m}}(X_{k-1},\overline{X}_{k-1}))+E^{Q_{n,x}}(D_{n-\overline{m}}(p_{n-\overline{m},\omega}(X_{k-1},\cdot), \overline{p}_{n-\overline{m}}(X_{k-1},\cdot))).\end{multline}

The Kantorovich-Rubinstein theorem and Control \ref{Holder} are used to bound the inequality's final term.  Let $f:\mathbb{R}^d\rightarrow\mathbb{R}$ be a function satisfying $\abs{f(x)-f(y)}\leq d_{n-\overline{m}}(x,y).$  Then, choose a smooth cutoff function $0\leq\chi_{n-\overline{m}}\leq 1$ satisfying $$\chi_{n-\overline{m}}(x)=\left\{\begin{array}{ll} 1 & \textrm{on}\;\; \overline{B}_{\tilde{D}_{n-\overline{m}}}, \\ 0 & \textrm{on}\;\;\mathbb{R}^d\setminus B_{2\tilde{D}_{n-\overline{m}}},\end{array}\right.$$ and for which $\abs{\chi_{n-\overline{m}}}_{n-\overline{m}}\leq 3$.

Fix $y\in[-L_{n+2}^2, L_{n+2}^2]^d$ and define $\tilde{f}_y(z)=f(z)-f(y)$.  Then, recalling the notation from (\ref{prob_operators}) and (\ref{prob_difference}), \begin{equation}\label{couple_main_30}\abs{S_{n-\overline{m}}f(y)}=\abs{S_{n-\overline{m}}\tilde{f}_y(y)}\leq\abs{S_{n-\overline{m}}(\chi_{n-\overline{m}}\tilde{f}_y)(y)}+\abs{S_{n-\overline{m}}((1-\chi_{n-\overline{m}})\tilde{f}_y)(y)}.\end{equation}  Since, for $C>0$ independent of $n$, $$\abs{\chi_{n-\overline{m}}\tilde{f}_y}_{n-\overline{m}}\leq C\tilde{\kappa}_{n-\overline{m}},$$ Control \ref{Holder}, which is satisfied owing to the assumptions $y\in[-L_{n+2}^2, L_{n+2}^2]^d$ and $\omega\in A_n$, implies \begin{equation}\label{couple_main_4}\abs{S_n(\chi_{n-\overline{m}}\tilde{f}_y)(y)}\leq CL_{n-\overline{m}}^{-\delta}\tilde{\kappa}_{n-\overline{m}}.\end{equation}  The second term is bounded using Control \ref{localization} since $\omega\in A_n$ and $y\in[-L_{n+2}^2, L_{n+2}^2]^d$.  First, by the triangle inequality,  $$\abs{S_{n-\overline{m}}((1-\chi_{n-\overline{m}})\tilde{f}_y)(y)}\leq\abs{R_{n-\overline{m}}((1-\chi_{n-\overline{m}})\tilde{f}_y)(y)}+\abs{\overline{R}_{n-\overline{m}}((1-\chi_{n-\overline{m}})\tilde{f}_y)(y)}.$$  Following an integration by parts and bounding $$P_{y,\omega}(\abs{X_{L_{n-\overline{m}}^2}}\geq r)\leq P_{y,\omega}(X^*_{L_{n-\overline{m}}^2}\geq r),$$ Control \ref{localization} yields the estimate, for $C>0$ independent of $n$, \begin{multline}\label{couple_main_5}\abs{R_{n-\overline{m}}((1-\chi_{n-\overline{m}})\tilde{f}_y)(y)}\leq -C \int_{\tilde{D}_{n-\overline{m}}}^\infty (\frac{r}{L_{n-\overline{m}}})^\beta \frac{d}{dr}P_{y,\omega}(\abs{X_{L_{n-\overline{m}}^2}}\geq r)\;dr\leq \\ C\tilde{\kappa}_{n-\overline{m}}^\beta\exp(-\tilde{\kappa}_{n-\overline{m}}).\end{multline}  Then, using the explicit representation of the Green's function and the control of $\alpha_{n-\overline{m}}$ present in Theorem \ref{effectivediffusivity}, for $C,c>0$ independent of $n$, \begin{equation}\label{couple_main_6}\abs{\overline{R}_{n-\overline{m}}((1-\chi_{n-\overline{m}})\tilde{f}_y)(y)}\leq \int_{\mathbb{R}^d\setminus\overline{B}_{\tilde{D}_{n-\overline{m}}}}\overline{p}_{n-\overline{m}}(y,z)(\frac{\abs{z-y}}{L_{n-\overline{m}}})^\beta\;dz\leq C\exp(-c\tilde{\kappa}_{n-\overline{m}}^2).\end{equation}  Since, in view of (\ref{L}), (\ref{kappa}) and (\ref{D}) there exists $C>0$ such that, for all $n\geq \overline{m}$, $$\exp(-c\tilde{\kappa}_{n-\overline{m}}^2)\leq C\tilde{\kappa}_{n-\overline{m}}^\beta\exp(-\tilde{\kappa}_{n-\overline{m}})\leq CL_{n-\overline{m}}^{-\delta}\tilde{\kappa}_{n-\overline{m}},$$ the combination (\ref{couple_main_30}), (\ref{couple_main_4}), (\ref{couple_main_5}) and (\ref{couple_main_6}) yields, for $C>0$ independent of $n$, \begin{equation}\label{couple_main_7}\abs{S_{n-\overline{m}}f(y)}\leq C\tilde{\kappa}_{n-\overline{m}}L_{n-\overline{m}}^{-\delta}.\end{equation}

If $y\notin[-L_{n+2}^2, L_{n+2}^2]^d$, again defining $\tilde{f}_y(z)=f(z)-f(y)$, \begin{equation}\label{couple_main_9}\abs{S_{n-\overline{m}}f(y)}=\abs{S_{n-\overline{m}}\tilde{f}_y(y)}\leq \abs{R_{n-\overline{m}}\tilde{f}_y(y)}+\abs{\overline{R}_{n-\overline{m}}\tilde{f}_y(y)}.\end{equation}  To bound the first term, recall that, almost surely with respect to $P_{y,\omega}$, paths $X_s\in\C([0,\infty);\mathbb{R}^d)$ satisfy the stochastic differential equation $$\left\{\begin{array}{l} dX_s=b(X_s,\omega)dt+\sigma(X_s,\omega)dB_s, \\ X_0=y,\end{array}\right.$$  for $A(x,\omega)=\sigma(x,\omega)\sigma(x,\omega)^t$ and for $B_t$ some standard Brownian motion under $P_{x,\omega}$ with respect to the canonical right-continuous filtration on $\C([0,\infty);\mathbb{R}^d)$.  Therefore, using the exponential inequality for Martingales, see Revuz and Yor \cite[Chapter~2, Proposition~1.8]{RY}, (\ref{bounded}) and (\ref{Lipschitz}), for every $R\geq 0$, for $C>0$ independent of $R$, $y$ and $\omega$, $$P_{y,\omega}(X^*_{L_{n-\overline{m}}^2}\geq R+CL_{n-\overline{m}}^2)\leq \exp(-\frac{R^2}{CL_{n-\overline{m}}^2}).$$  Choosing $R=C\tilde{\kappa}_{n-\overline{m}}L_{n-\overline{m}}$, it follows that, for $C,c>0$ independent of $n\geq\overline{m}$, $$P_{y,\omega}(X^*_{L_{n-\overline{m}}^2}\geq C\tilde{\kappa}_{n-\overline{m}}L_{n-\overline{m}}^2)\leq \exp(-c\tilde{\kappa}_{n-\overline{m}}^2).$$  Then, form the decomposition \begin{multline*} \abs{R_{n-\overline{m}}\tilde{f}_y(y)}\leq\abs{E_{y,\omega}(\tilde{f}_y(X_{L_{n-\overline{m}}^2}), X^*_{L_{n-\overline{m}}^2}\leq C\tilde{\kappa}_{n-\overline{m}}L_{n-\overline{m}}^2)} \\ +\abs{E_{y,\omega}(\tilde{f}_y(X_{L_{n-\overline{m}}^2}), X^*_{L_{n-\overline{m}}^2}> C\tilde{\kappa}_{n-\overline{m}}L_{n-\overline{m}}^2)}.\end{multline*}  The second term of this inequality is bounded using (\ref{couple_main_5}).  The first term is bounded brutally, using the fact that $\tilde{f}(z)\leq d_{n-\overline{m}}(z,y)$, which yields, for $C>0$ independent of $n$, $$\abs{P_{y,\omega}(\tilde{f}_y(X_{L_{n-\overline{m}}^2}), X^*_{L_{n-\overline{m}}^2}\leq C\tilde{\kappa}_{n-\overline{m}}L_{n-\overline{m}}^2)}\leq C\tilde{\kappa}_{n-\overline{m}}^\beta L_{n-\overline{m}}^\beta.$$  Therefore, for $C>0$ independent of $n$, \begin{equation}\label{couple_main_8}\abs{R_{n-\overline{m}}\tilde{f}_y(y)}\leq C\tilde{\kappa}_{n-\overline{m}}^\beta L_{n-\overline{m}}^\beta.\end{equation}  The second term of (\ref{couple_main_9}) is bounded using the explicit representation of the heat kernel and Theorem \ref{effectivediffusivity}.  For $C>0$ independent of $n$, since $\tilde{f}_y(z)\leq d_{n-\overline{m}}(z,y)$, \begin{equation}\label{couple_main_10}\abs{\overline{R}_{n-\overline{m}}\tilde{f}_y(y)}\leq \abs{\int_{\mathbb{R}^d}\overline{p}_{n-\overline{m}}(y,z)(\frac{\abs{z-y}}{L_{n-\overline{m}}})^\beta\;dz}\leq C.\end{equation}  Therefore, in view of (\ref{couple_main_9}), (\ref{couple_main_8}) and (\ref{couple_main_10}), for $C>0$ independent of $n$, \begin{equation}\label{couple_main_11} \abs{S_nf(y)}\leq C(1+\tilde{\kappa}_{n-\overline{m}}^\beta L_{n-\overline{m}}^{\beta})\leq C\tilde{\kappa}_{n-\overline{m}}^\beta L_{n-\overline{m}}^{\beta}.\end{equation}

Returning to (\ref{couple_main_3}), decompose the second term as \begin{multline}\label{couple_main_12}E^{Q_{n,x}}(D_{n-\overline{m}}(p_{n-\overline{m},\omega}(X_{k-1},\cdot), \overline{p}_{n-\overline{m}}(X_{k-1},\cdot)))= \\ E^{Q_{n,x}}(D_{n-\overline{m}}(p_{n-\overline{m},\omega}(X_{k-1},\cdot), \overline{p}_{n-\overline{m}}(X_{k-1},\cdot)), X_{k-1}\in[L_{n+2}^2, L_{n+2}^2]^d) \\ +E^{Q_{n,x}}(D_{n-\overline{m}}(p_{n-\overline{m},\omega}(X_{k-1},\cdot), \overline{p}_{n-\overline{m}}(X_{k-1},\cdot)), X_{k-1}\notin[L_{n+2}^2, L_{n+2}^2]^d).\end{multline}  Since $f$ satisfying $\abs{f(x)-f(y)}\leq d_{n-\overline{m}}(x,y)$ appearing in (\ref{couple_main_7}) and (\ref{couple_main_11}) was arbitrary, the Kantorovich-Rubinstein theorem in (\ref{couple_kr}) with (\ref{couple_main_7}) imply that \begin{equation}\label{couple_main_14}E^{Q_{n,x}}(D_{n-\overline{m}}(p_{n-\overline{m},\omega}(X_{k-1},\cdot), \overline{p}_{n-\overline{m}}(X_{k-1},\cdot)), X_{k-1}\in[L_{n+2}^2, L_{n+2}^2]^d)\leq C\tilde{\kappa}_{n-\overline{m}}L_{n-\overline{m}}^{-\delta}.\end{equation}  Then, again the Kantorovich-Rubinstein theorem from (\ref{couple_kr}) together with (\ref{couple_main_11}), using  $x\in[-\frac{1}{2}L_{n+2}^2, \frac{1}{2}L_{n+2}^2]^d$ and $0\leq k\leq 2(\frac{L_{n+2}}{L_{n-\overline{m}}})^d$, the second term is bounded by \begin{multline}\label{couple_main_15}E^{Q_{n,x}}(D_{n-\overline{m}}(p_{n-\overline{m},\omega}(X_{k-1},\cdot), \overline{p}_{n-\overline{m}}(X_{k-1},\cdot)), X_{k-1}\notin[L_{n+2}^2, L_{n+2}^2]^d) \\ \leq C\tilde{\kappa}_{n-\overline{m}}^\beta L_{n-\overline{m}}^{\beta}P_{x,\omega}({X}_{2L_{n+2}^2}^*\geq \frac{1}{2}L_{n+2}^2).\end{multline}  To conclude, since $x\in[-\frac{1}{2}L_{n+2}^2, \frac{1}{2}L_{n+2}^2]^d$ and $\omega\in A_n$, Control \ref{localization} implies that, for $C>0$ independent of $n$, \begin{equation}\label{couple_main_2000}  P_{x,\omega}({X}_{2L_{n+2}^2}^*\geq \frac{1}{2}L_{n+2}^2)\leq 2\exp(-\frac{L_{n+2}^2}{D_{n+2}})\leq C\exp(-L_{n+1}).\end{equation}  Therefore, in view of (\ref{couple_main_12}), (\ref{couple_main_14}), (\ref{couple_main_15}) and (\ref{couple_main_2000}), for $C>0$ and $c>0$ independent of $n$, \begin{multline}\label{couple_main_17} E^{Q_{n,x}}(D_{n-\overline{m}}(p_{n-\overline{m},\omega}(X_{k-1},\cdot), \overline{p}_{n-\overline{m}}(X_{k-1},\cdot))) \\ \leq \tilde{\kappa}_{n-\overline{m}}L_{n-\overline{m}}^{-\delta}+C\tilde{\kappa}_{n-\overline{m}}^\beta L_{n-\overline{m}}^{\beta}\exp(-L_{n+1})\leq C\tilde{\kappa}_{n-\overline{m}}L_{n-\overline{m}}^{-\delta}.\end{multline}

It follows from (\ref{couple_main_3}) and (\ref{couple_main_12}) that, for every $0\leq k\leq 2(\frac{L_{n+2}}{L_{n-\overline{m}}})^2$, for $C>0$ independent of $n$, $$E^{Q_{n,x}}(d_{n-\overline{m}}(X_k,\overline{X}_k))\leq E^{Q_{n,x}}(d_{n-\overline{m}}(X_{k-1},\overline{X}_{k-1}))+C\tilde{\kappa}_{n-\overline{m}}L_{n-\overline{m}}^{-\delta}.$$  And, since the initial distribution of the Markov chain $(X_0,\overline{X}_0)=(x,x)$ with probability one under $Q_{n,x}$, iterating this inequality yields, for every $0\leq k\leq 2(\frac{L_{n+2}}{L_{n-\overline{m}}})^2$,  $$E^{Q_{n,x}}(d_{n-\overline{m}}(X_k,\overline{X}_k))\leq Ck\tilde{\kappa}_{n-\overline{m}}L_{n-\overline{m}}^{-\delta}.$$  Chebyshev's inequality and the definition of $d_{n-\overline{m}}$ then imply that, for every $0\leq k\leq 2(\frac{L_{n+2}}{L_{n-\overline{m}}})^2$, for every $\gamma>0$, $$(\frac{\gamma}{L_{n-\overline{m}}})^\beta Q_{n,x}(\abs{X_k-\overline{X}_k}\geq \gamma)\leq Ck\tilde{\kappa}_{n-\overline{m}}L_{n-\overline{m}}^{-\delta},$$ and, therefore, \begin{multline*}Q_{n,x}(\abs{X_k-\overline{X}_k}\geq \gamma\;|\;\textrm{for some}\;0\leq k\leq 2(\frac{L_{n+2}}{L_{n-\overline{m}}})^2)\leq \\ C(\frac{L_{n-\overline{m}}}{\gamma})^\beta(1+\ldots+2(\frac{L_{n+2}}{L_{n-\overline{m}}})^2)\tilde{\kappa}_{n-\overline{m}}L_{n-\overline{m}}^{-\delta}.\end{multline*}  Estimating the sum by the elementary inequality $1+\ldots+m\leq m^2$, $$Q_{n,x}(\abs{X_k-\overline{X}_k}\geq \gamma\;|\;\textrm{for some}\;0\leq k\leq 2(\frac{L_{n+2}}{L_{n-\overline{m}}})^2)\leq C(\frac{L_{n-\overline{m}}}{\gamma})^\beta(\frac{L_{n+2}}{L_{n-\overline{m}}})^4\tilde{\kappa}_{n-\overline{m}}L_{n-\overline{m}}^{-\delta},$$ which, since $n\geq\overline{m}$, $\omega\in A_n$ and $x\in[-\frac{1}{2}L_{n+2}^2, \frac{1}{2}L_{n+2}^2]^d$ were arbitrary, completes the argument. \end{proof}

The section concludes with a straightforward corollary of Proposition \ref{couple_main}.  Since the definition of $\overline{m}$ in (\ref{prob_m}) implies $$(1+a)^{\overline{m+2}}-1\leq \frac{(1+a)^4}{2-(1+a)^2}-1\leq \frac{8a}{2}=4a,$$ it follows from the definition of $L_n$ in (\ref{L}) that, for $C>0$ independent of $n$, $$(\frac{L_{n+2}}{L_{n-\overline{m}}})^4\leq C L_{n-\overline{m}}^{16a}.$$  The following corollary then follows immediately by taking $\gamma=L_{n-\overline{m}}$ in Proposition \ref{couple_main}.  Observe that (\ref{Holderexponent}) and (\ref{delta}) imply the exponent $$16a-\delta<0$$ is negative.

\begin{cor}\label{couple_cor}  Assume (\ref{steady}) and (\ref{constants}).  For every $n\geq \overline{m}$, $\omega\in A_n$ and $x\in[-\frac{1}{2}L_{n+2}^2, \frac{1}{2}L_{n+2}^2]^d$, for $C>0$ independent of $n$, $$Q_{n,x}(\abs{X_k-\overline{X}_k}\geq L_{n-\overline{m}}\;|\;\textrm{for some}\;0\leq k\leq 2(\frac{L_{n+2}}{L_{n-\overline{m}}})^2)\leq C\tilde{\kappa}_{n-\overline{m}}L_{n-\overline{m}}^{16a-\delta}.$$\end{cor}

\section{Estimates for the Exit Time of Brownian Motion}\label{Brownian_exit}

In this section estimates are obtained, in expectation and near the boundary of the domain, for the exit time of a Brownian motion.  The role of the exterior ball condition comes in the proof of these estimates, which states that there exists (a now fixed) $r_0>0$ such that, for every $x\in\partial U$, there exists $x^*\in\mathbb{R}^d$ satisfying \begin{equation}\label{disc_ball} \overline{B}_{r_0}(x^*)\cap\overline{U}=\left\{x\right\}.\end{equation}  Furthermore, define, for each $\delta>0$, the inflated domain \begin{equation}\label{disc_inflate} U_{\delta}=\left\{\;x\in\mathbb{R}^d\;|\;d(x,U)< \delta\;\right\},\end{equation} and notice, as a consequence of (\ref{disc_ball}), for every $0<\delta<r_0$, \begin{equation}\label{disc_inflate_ball} U_{\delta}\;\;\textrm{satisfies the exterior ball condition with radius}\;(r_0-\delta).\end{equation}  Essentially, it will be necessary to understand, in expectation, the exit time of Brownian motion from the sets $U_\delta$ and $U$, as $\delta\rightarrow 0$, at points within distance $\delta$ from the boundary.

The first step is to consider the exit time of Brownian motion from the annular domains centered at the origin and defined, for each pair of radii $0<r_1<r_2<\infty$, by $$A_{r_1, r_2}=B_{r_2}\setminus\overline{B}_{r_1}.$$  For each pair $(r_1,r_2)$ let $\tau_{r_1,r_2}$ denote the $\C([0,\infty);\mathbb{R}^d)$ exit time $$\tau_{r_1,r_2}=\inf\left\{\;t\geq 0\;|\;X_t\notin A_{r_1,r_2}\;\right\},$$ and recall that, in expectation and with respect to the Weiner measure $W^{n}_x$ defining Brownian motion with variance $\alpha_n$ beginning from $x$, the function $$u^n_{r_1,r_2}(x)=E^{W^n_x}(\tau_{r_1,r_2})\;\;\textrm{on}\;\;\overline{A}_{r_1,r_2}$$ satisfies the equation \begin{equation}\label{disc_annulus_exit}\left\{\begin{array}{ll} 1+\frac{\alpha_n}{2}\Delta u^n_{r_1,r_2}=0 & \textrm{on}\;\; A_{r_1,r_2}, \\ u^n_{r_1,r_2}=0 & \textrm{on}\;\;\partial A_{r_1,r_2}.\end{array}\right.\end{equation}  See, for instance, \cite[Exercise~9.12]{Oksendal}.  The following proposition obtains an upper bound for these solutions most effective in a neighborhood of $\partial B_{r_1}$.  The estimate necessarily depends upon the pair $(r_1,r_2)$, which in the application to follow will be fixed independently of $n\geq 0$, but does not otherwise depend upon $n\geq 0$.

\begin{prop}\label{disc_annulus}  Assume (\ref{steady}) and (\ref{constants}).  For each pair of radii $0<r_1<r_2<\infty$, for each $n\geq 0$, there exists $C=C(r_1,r_2)>0$ such that $$u^n_{r_1,r_2}(x)\leq C\d(x,\partial B_{r_1})\;\;\textrm{on}\;\;\overline{A}_{r_1,r_2}.$$\end{prop}

\begin{proof}  Fix $0<r_1<r_2<\infty$ and $n\geq 0$.  The solution $u^n_{r_1,r_2}$ of (\ref{disc_annulus_exit}) admits the explicit radial description, since $d\geq 3$ owing to (\ref{dimension}), writing $r=\abs{x}$, $$u(x)=u(r)=c_1(r_1,r_2)+c_2(r_1,r_2)r^{2-d}-\frac{r^2}{2d\alpha_n}\;\;\textrm{on}\;\;A_{r_1,r_2},$$ where $$c_1(r_1,r_2)=\frac{1}{2d\alpha_n}\cdot \frac{r_1^2r_2^{2-d}-r_2^2r_1^{2-d}}{r_2^{2-d}-r_1^{2-d}}\;\;\textrm{and}\;\;c_2(r_1,r_2)=\frac{1}{2d\alpha_n}\cdot\frac{r_2^2-r_1^2}{r_2^{2-d}-r_1^{2-d}}.$$  Therefore, after performing a Taylor expansion in $r$ about $r_1$, using the fact that $u(r_1)=0$, for each $x\in A_{r_1,r_2}$, $$u^n_{r_1,r_2}(x)=u^n_{r_1,r_2}(r)=c_2(2-d)r_1^{1-d}(r-r_1)+c_2(2-d)(1-d)\int_{r_1}^rs^{-d}(r-s)\;ds-(\frac{2r_1r+r^2}{2d\alpha_n}).$$  Because the integrand is bounded by $r_1^{-d}(r_2-r_1)$, a brutal but sufficient estimate, and because the final term is negative, the uniform control in $n\geq 0$ of $\alpha_n$ provided by Theorem \ref{effectivediffusivity} implies that there exits $C=C(r_1,r_2)>0$ satisfying $$u^n_{r_1,r_2}(x)=u^n_{r_1,r_2}(r)\leq C(r-r_1)=C\d(x,\partial B_{r_1})\;\;\textrm{on}\;\;\overline{A}_{r_1,r_2},$$ which completes the argument.\end{proof}

Passing from the annular regions $A_{r_1,r_2}$ to the domain $U$ and its inflations $U_\delta$, for $\delta>0$ small, is now straightforward.  Define, for each $x\in\mathbb{R}^d$ and pair of radii $(r_1,r_2)$, the translate $$A_{r_1,r_2}(x)=x+A_{r_1,r_2}=B_{r_2}(x)\setminus\overline{B}_{r_1}(x),$$ and, for each $\delta>0$, the $\C([0,\infty);\mathbb{R}^d)$ exit times \begin{equation}\label{disc_stopping} \tau=\inf\left\{\;t\geq 0\;|\;X_t\notin U\;\right\}\;\;\textrm{and}\;\;\tau^\delta=\inf\left\{\; t\geq 0\;|\;X_t\notin U_{\delta}\;\right\}.\end{equation}  The following corollary controls the expectation of $\tau$ and $\tau^\delta$ in an approximately $\delta$-neighborhood of the respective boundaries.  Recall the radius $r_0$ in (\ref{disc_ball}) quantifying the exterior ball condition.

\begin{cor}\label{disc_u} Assume (\ref{steady}) and (\ref{constants}).  For every $0<\delta<\frac{r_0}{2}$, for every $n\geq 0$, for $C>0$ independent of $n$ and $\delta$, $$\sup_{\d(x,\partial U)\leq \delta} E^{W^n_{x}}(\tau)\leq C\delta\;\;\textrm{and}\;\;\sup_{\d(x,\partial U_\delta)\leq 2\delta}E^{W^n_x}(\tau^\delta)<C\delta.$$\end{cor}

\begin{proof}  For each $0<\delta<\frac{r_0}{2}$ it follows from (\ref{disc_inflate_ball}) that $U_{\delta}$ satisfies the exterior ball condition with radius $r_0-\delta$.  Fix $r_2>\frac{r_0}{2}$ such that, whenever $x\in\partial U_\delta$ and $x^*\in\mathbb{R}^d$ satisfy $$\overline{B}_{r_0-\delta}(x^*)\cap\overline{U}_\delta=\left\{x\right\},\;\;\textrm{it follows that}\;\;\overline{U}_\delta\subset B_{r_2}(x^*).$$  The existence of $r_2$ chosen uniformly for $0<\delta<\frac{r_0}{2}$ is guaranteed by the boundedness of $U$ assumed in (\ref{domain_bounded}).  Observe that, for each $x\in U$, since the stopping time $\tau_{\delta}\geq \tau$ almost-surely with respect to $W^n_x$, the first statement is subsumed by the second, which will be shown henceforth.

Fix $0<\delta<\frac{r_0}{2}$ and $n\geq 0$.  In Proposition \ref{disc_annulus}, choose $r_1=r_0-\delta$ and observe that the smallness of $\delta$ guarantees $\frac{r_0}{2}\leq r_1\leq r_0$, and the choice of pair $(r_1,r_2)$ ensures, for every $x\in\partial U_\delta$ and $x^*\in\mathbb{R}^d$ satisfying \begin{equation}\label{disc_u_1}\overline{B}_{r_0-\delta}(x^*)\cap\overline{U}_\delta=\left\{x\right\},\;\textrm{the containment}\;\;\overline{U}_\delta\subset \overline{A}_{r_1,r_2}(x^*).\end{equation}

Fix $x\in U_\delta$ satisfying $\d(x, \partial U_\delta)\leq 2\delta$ and, owing to the compactness of $\partial U_\delta$, choose $\overline{x}\in\partial U_\delta$ satisfying $\abs{x-\overline{x}}=\d(x, \partial U_\delta)$.  Let $\overline{x}^*$ satisfy (\ref{disc_u_1}) corresponding to $\overline{x}$, and let $u^{n,\overline{x}}_{r_1,r_2}$ satisfy the equation \begin{equation}\label{disc_u_3}\left\{\begin{array}{ll} 1+\frac{\alpha_n}{2}\Delta u^{n,\overline{x}}_{r_1,r_2}=0 & \textrm{on}\;\;A_{r_1,r_2}(\overline{x}^*), \\ u^{n,\overline{x}}_{r_1,r_2}=0 & \textrm{on}\;\;\partial A_{r_1,r_2}(\overline{x}^*).\end{array}\right.\end{equation}  Due to the translational invariance of the equation, Proposition \ref{disc_annulus} implies that, since $r_1$ and $r_2$ can be chosen to be bounded uniformly above and away from zero independently of $0<\delta<\frac{r_0}{2}$ and $n\geq 0$, for $C>0$ independent of $n$ and $\delta$, \begin{equation}\label{disc_u_2}u^{n,\overline{x}}_{r_1,r_2}(x)\leq C\d(x,\partial A_{r_1,r_2}(\overline{x}^*))\leq C\abs{x-\overline{x}}=C\d(x, \partial U_\delta)\leq C\delta.\end{equation}

In order to conclude, the expectation of the exit time $E^{W^n_x}(\tau^\delta)$ on $\overline{U}_\delta$ is the solution of the equation $$\left\{\begin{array}{ll} 1+\frac{\alpha_n}{2}\Delta E^{W^n_x}(\tau^\delta)=0 & \textrm{on}\;\; U_\delta,\\ E^{W^n_x}(\tau^\delta)=0 & \textrm{on}\;\;\partial U_\delta,\end{array}\right.$$ see \cite[Exercise~9.12]{Oksendal}, and (\ref{disc_u_1}) implies that the solution $u^{n,\overline{x}}_{r_1,r_2}\geq 0$ of (\ref{disc_u_3}) satisfies $$\left\{\begin{array}{ll} 1+\frac{\alpha_n}{2}\Delta u^{n,\overline{x}}_{r_1,r_2}(\tau_\delta)=0 & \textrm{on}\;\; U_\delta,\\ u^{n,\overline{x}}_{r_1,r_2}(\tau_\delta)\geq0 & \textrm{on}\;\;\partial U_\delta.\end{array}\right.$$  Therefore, by comparison principle and (\ref{disc_u_2}), for $C>0$ independent of $n$ and $\delta$, $$E^{W^n_x}(\tau^\delta)\leq u^{n,\overline{x}}_{r_1,r_2}(x)\leq C\delta,$$ which, because $x\in U_\delta$ satisfying $\d(x,\partial U_\delta)\leq 2\delta$, $0<\delta<\frac{r_0}{2}$ and $n\geq 0$ were arbitrary, completes the proof.\end{proof}

Corollary \ref{disc_u} is also sufficient to estimate, for each $\epsilon>0$, the exit times of Brownian motion from the rescaled domains $U/\epsilon$.  Write, for each $\epsilon>0$, the $\C([0,\infty);\mathbb{R}^d)$ exit time \begin{equation}\label{disc_epsilon_exit}\tau^\epsilon=\inf\left\{\;t\geq 0\;|\;X_t\notin U/\epsilon\right\}=\inf\left\{\;t\geq 0\;|\;\epsilon X_t\notin U\right\}.\end{equation}  The corresponding expectation $E^{W^n_x}(\tau^\epsilon)$ then satisfies $$\left\{\begin{array}{ll} 1+\frac{\alpha_n}{2}\Delta E^{W^n_x}(\tau^\epsilon)=0 & \textrm{on}\;\;U/\epsilon, \\ E^{W^n_x}(\tau^\epsilon)=0 & \textrm{on}\;\;\partial U/\epsilon,\end{array}\right.$$ see \cite[Exercise~9.12]{Oksendal}, and can be obtained by a rescaling of $E^{W^n_x}(\tau)$ from Corollary \ref{disc_u}.  Indeed, $$E^{W^n_x}(\tau^\epsilon)=\epsilon^{-2}E^{W^n_{\epsilon x}}(\tau)\;\;\textrm{on}\;\;\overline{U}/\epsilon.$$  Similarly, if for each $\epsilon>0$ and $\delta>0$, $$\tau^{\epsilon,\delta}=\inf\left\{\;t\geq 0\;|\;X_t\notin U_\delta/\epsilon\;\right\}=\inf\left\{\;t\geq 0\;|\;\epsilon X_t\notin U_\delta\;\right\},$$ then, for $\tau^\delta$ the exit time from $U_\delta$ defined in (\ref{disc_stopping}), $$E^{W^n_x}(\tau^{\epsilon,\delta})=\epsilon^{-2}E^{W^n_{\epsilon x}}(\tau^\delta)\;\;\textrm{on}\;\;\overline{U}_\delta/\epsilon.$$  The following statement is then an immediate consequence Corollary \ref{disc_u} and the previous two equalities.

\begin{cor}\label{disc_u_scale}  Assume (\ref{steady}) and (\ref{constants}).  For every $\epsilon>0$, $0<\delta<\frac{r_0}{2\epsilon}$ and $n\geq 0$, for $C>0$ independent of $\epsilon$, $\delta$ and $n$, $$\sup_{\d(x,\partial U/\epsilon)\leq \delta} E^{W^n_{x}}(\tau^\epsilon)\leq C\epsilon^{-1}\delta\;\;\textrm{and}\;\;\sup_{\d(x,\partial U_{\delta}/\epsilon)\leq 2\delta}E^{W^n_x}(\tau^{\epsilon, \delta})<C\epsilon^{-1}\delta.$$\end{cor}

\section{The Discrete Approximation and Proof of Theorem \ref{intro_main}}\label{main}

The purpose of this section is to complete the almost-sure characterization, as $\epsilon\rightarrow 0$, of solutions \begin{equation}\label{end_eq_100}\left\{\begin{array}{ll} \frac{1}{2}\tr(A(\frac{x}{\epsilon},\omega)D^2u^\epsilon)+\frac{1}{\epsilon}b(\frac{x}{\epsilon},\omega)\cdot Du^\epsilon=0 & \textrm{on}\;\;U, \\ u^\epsilon=f(x) & \textrm{on}\;\;\partial U,\end{array}\right.\end{equation} which, for the $\C([0,\infty);\mathbb{R}^d)$ stopping time \begin{equation}\label{end_tau}\tau^\epsilon=\inf\left\{\;t\geq 0\;|\;\epsilon X_t\notin U\;\right\}=\inf\left\{\;t\geq 0\;|\;X_t\notin U/\epsilon\;\right\},\end{equation} and probability measure and expectation $P_{x,\omega}$ and $E_{x,\omega}$ on $\C([0,\infty);\mathbb{R}^d)$ associated to the unscaled generator \begin{equation}\label{end_gen}\frac{1}{2}\sum_{i,j=1}^da_{ij}(x,\omega)\frac{\partial^2}{\partial x_i\partial x_j}+\sum_{i=1}^db_i(x,\omega)\frac{\partial}{\partial x_i},\end{equation} admit the representation $$u^\epsilon(x)=E_{\frac{x}{\epsilon},\omega}(f(\epsilon X_{\tau^\epsilon}))\;\;\textrm{on}\;\;\overline{U}.$$  See (\ref{intro_rep}), (\ref{intro_res}) and \cite[Exercise~9.12]{Oksendal}.

The strategy will be, for scales $\epsilon$ satisfying $L_n\leq \frac{1}{\epsilon}<L_{n+1}$, to approximate the continuous process $X_t$ by the discrete process constructed in Proposition \ref{couple_main} corresponding to time steps of order $L_{n-\overline{m}}^2$.  The choice in (\ref{prob_m}) of the integer $\overline{m}$ satisfying $$\overline{m}>1-\frac{\log(1-2a-a^2)}{\log(1+a)}$$ guarantees, in view of the definitions of $L_n$ in (\ref{L}) and $\tilde{D}_n$ in (\ref{D}), that there exists $\zeta>0$ such that, for $C>0$ independent of $n$, \begin{equation}\label{end_zeta_1}L_{n+1}\tilde{D}_{n-\overline{m}}\leq CL_{n-1}^{2-\zeta}.\end{equation}  Furthermore, in order to simplify some statements to follow, use the definition of $L_n$ in (\ref{L}) and $\tilde{\kappa}_n$ in (\ref{kappa}) to fix $\zeta>0$ sufficiently small so that, for all $n\geq\overline{m}$, \begin{equation}\label{end_zeta_2}\tilde{\kappa}_{n-\overline{m}}L_{n-\overline{m}}^{16a-\delta}\leq L_{n-1}^{-\zeta}\;\;\textrm{and}\;\;L_n^{-da}\leq L_{n-1}^{-\zeta}.\end{equation}  Henceforth, \begin{equation}\label{end_zeta} \textrm{fix}\;\;\zeta>0\;\;\textrm{satisfying (\ref{end_zeta_1}) and (\ref{end_zeta_2}).}\end{equation}

Introduce, for each $\epsilon>0$ and $n\geq \overline{m}$, the discrete $\C([0,\infty);\mathbb{R}^d)$ stopping times \begin{equation}\label{end_tau1} \tau^{\epsilon,n}_1=\inf\left\{\;kL_{n-\overline{m}}^2\geq 0\;|\;\d(X_{kL_{n-\overline{m}}^2},(U/\epsilon)^c)\leq \tilde{D}_{n-\overline{m}}\;\right\},\end{equation} which quantify the first time in the discrete sequence $\left\{kL_{n-\overline{m}}^2\right\}_{k\geq 0}$ that $X_{kL_{n-\overline{m}}^2}$ resides in the $\tilde{D}_{n-\overline{m}}$ neighborhood of the compliment of $U/\epsilon$.  It is not true that $\tau^{\epsilon,n}\leq \tau^\epsilon$ for every path $X_t$, however, for scales $L_n\leq \frac{1}{\epsilon}<L_{n+1}$, the failure of this inequality will be controlled in probability by the exponential estimate appearing Control \ref{localization}.

Furthermore, define, for each $\epsilon>0$ and $n\geq \overline{m}$, the discrete $\C([0,\infty);\mathbb{R}^d)$ stopping times \begin{equation}\label{end_tau2} \tau^{\epsilon,n}_2=\inf\left\{\;kL_{n-\overline{m}}^2\geq 0\;|\;\d(X_{kL_{n-\overline{m}}^2},(U/\epsilon))\geq \tilde{D}_{n-\overline{m}}\;\right\}.\end{equation}  These stopping times indicate the first earliest point in the discrete sequence $\left\{kL_{n-\overline{m}}^2\right\}_{k\geq 0}$ that $X_{kL_{n-\overline{m}}^2}$ resides outside the $\tilde{D}_{n-\overline{m}}$ neighborhood of $(U/\epsilon)$.

It follows from the definitions that $\tau^{\epsilon,n}_1\leq \tau^{\epsilon,n}_2$ and, on the event $\tau^{\epsilon,n}_1\leq \tau^\epsilon$, it is necessarily the case that $\tau^{\epsilon,n}_1\leq \tau^\epsilon\leq \tau^{\epsilon,n}_2.$  The estimates obtained in expectation for the exit time of Brownian motion appearing in Corollary \ref{disc_u_scale} will allow for an effective estimate of $\tau^{\epsilon,n}_2$ near the boundary of $U/\epsilon$ for scales $\epsilon$ and $n$ satisfying $L_n\leq\frac{1}{\epsilon}<L_{n+1}$.  These bounds, together with the coupling constructed in Proposition \ref{couple_main}, then yield an upper bound for the probability $$P_{x,\omega}(\tau^\epsilon-\tau_1^{\epsilon,n}\geq L_{n-1}^2)\;\;\textrm{for}\;\;x\in\overline{U}/\epsilon.$$

This estimate, in conjunction with the exponential controls established by Control \ref{localization}, effectively provides a barrier for equation (\ref{end_eq_100}) near of boundary $\partial U$ of a quality that, for general such equations, is impossible to obtain.  And, therefore, shows that the discretely stopped process $X_{\tau^{\epsilon,n}_1}$ is a good proxy for $X_{\tau^\epsilon}$.  The proof of Theorem \ref{intro_main} then follows from the coupling established in Corollary \ref{disc_u_scale} and the upper bound for the exit time appearing in Proposition \ref{exit_main}.

\begin{prop}\label{end_Brownian}  Assume (\ref{steady}) and (\ref{constants}).  For all $n\geq 0$ sufficiently large, for every $\epsilon>0$ satisfying $L_n\leq\frac{1}{\epsilon}\leq L_{n+1}$, for $\zeta>0$ in (\ref{end_zeta}) and $C>0$ independent of $n$, $$\sup_{\d(x,(U/\epsilon)^c)\leq 2\tilde{D}_{n-\overline{m}}}W^{n-\overline{m}}_x(\tau_2^{\epsilon,n}\geq L_{n-1}^2)\leq CL_{n-1}^{-\zeta}.$$\end{prop}

\begin{proof}  Fix $n_1\geq 0$ such that, for every $n\geq n_1$, for $r_0$ the constant quantifying the exterior ball condition in (\ref{disc_ball}), \begin{equation}\label{end_Brownian_1} 2\tilde{D}_{n-\overline{m}}\leq \frac{r_0L_n}{2}.\end{equation}  This condition guarantees that, whenever $\d(x,(U/\epsilon)^c)\leq 2\tilde{D}_{n-\overline{m}}$, the conclusion of Proposition \ref{disc_u_scale} is satisfied.

Henceforth, fix $n\geq n_1$, $\epsilon>0$ satisfying $L_n\leq \frac{1}{\epsilon}<L_{n+1}$ and $x\in \mathbb{R}^d$ satisfying $d(x,(U/\epsilon)^c)\leq2\tilde{D}_{n-\overline{m}}$.  Recalling that $\tau^{\epsilon,\delta}$ denotes the exit time from the $\delta$-neighborhood of $(U/\epsilon)$, and choosing $\delta=2\tilde{D}_{n-\overline{m}}$, Proposition \ref{disc_u_scale} implies that, for $C>0$ independent of $n$, $$E^{W^n_x}(\tau^{\epsilon,2\tilde{D}_{n-\overline{m}}})\leq C(2\tilde{D}_{n-\overline{m}})\epsilon^{-1}< C\tilde{D}_{n-\overline{m}}L_{n+1}.$$  Therefore, for $\zeta>0$ defined in (\ref{end_zeta}), for $C>0$ independent of $n$, $$E^{W^n_x}(\tau^{\epsilon,2\tilde{D}_{n-\overline{m}}})\leq CL_{n-1}^{2-\zeta}.$$  Then, by Chebyshev's inequality, for $C>0$ independent of $n$, \begin{equation}\label{end_Brownian_2} W^n_x(\tau^{\epsilon,2\tilde{D}_{n-\overline{m}}}\geq \frac{1}{2}L_{n-1}^2)\leq CL_{n-1}^{-\zeta}.\end{equation}

In order to conclude, using the translational invariance of the heat kernel and the Markov property, and owing to exponential tail estimates for Brownian motion on scale $\tilde{D}_{n-\overline{m}}$, see \cite[Chapter~2, Proposition~1.8]{RY}, using the choice of constants (\ref{L}), (\ref{kappa}) and (\ref{D}), for $C>0$ and $c>0$ independent of $n$, \begin{equation}\label{end_Brownian_3}W^n_x(\tau^{\epsilon,2\tilde{D}_{n-\overline{m}}}+L_{n-\overline{m}}^2\leq \tau^{\epsilon,n}_2)\leq W_0^n(X^*_{L_{n-\overline{m}}^2}\geq \tilde{D}_{n-\overline{m}})\leq C\exp(-c\tilde{\kappa}_{n-\overline{m}}^2).\end{equation}  And, since the choice of constants (\ref{L}), (\ref{kappa}) and (\ref{D}) guarantee the existence of $C>0$ independent of $n$ satisfying $\exp(-c\tilde{\kappa}_{n-\overline{m}}^2)\leq CL_{n-1}^{-\zeta},$ and since for $n\geq 0$ sufficiently large $L_{n-\overline{m}}^2<\frac{1}{2}L_{n-1}^2$, in combination (\ref{end_Brownian_2}) and (\ref{end_Brownian_3}) assert that $$W^n_x(\tau^{\epsilon,n}_2\geq L_{n-1}^2)\leq CL_{n-1}^{-\zeta},$$ which, since $x$ satisfying $d(x(U/\epsilon)^c)\leq2\tilde{D}_{n-\overline{m}}$, $\epsilon$ satisfying $L_n\leq \frac{1}{\epsilon}<L_{n+1}$ and $n\geq n_1$ were arbitrary, completes the argument.  \end{proof}

The following proposition relies upon the random subsets $A_n$ defined in (\ref{prob_event_1}).  Recall, for each $n\geq 0$, \begin{multline}\label{end_event}A_n=\left\{\;\omega\in\Omega\;|\;\omega\in B_m(x)\;\;\textrm{for all}\;\;x\in L_m\mathbb{Z}^d\cap[-L_{n+2}^2, L_{n+2}^2]^d\;\;\textrm{and}\right. \\ \left. \textrm{for all}\;\;n-\overline{m}\leq m\leq n+2.\right\},\end{multline} which implies that for every $\omega\in A_n$, $x\in[-L_{n+2}^2, L_{n+2}^2]^d$ and scale between $L_{n-\overline{m}}$ to $L_{n+2}$, the H\"older estimate from Control \ref{Holder} and the localization estimate from Control \ref{localization} are satisfied.  The remaining arguments require no further use of Control \ref{Holder}, since the coupling obtained in Corollary \ref{couple_cor} encodes already its purpose, but the localization estimate from Control \ref{localization} will be used, and which is now recalled.

\begin{con}\label{end_localization}  Fix $x\in\mathbb{R}^d$, $\omega\in\Omega$ and $n\geq 0$.  For each $v\geq D_n$, for all $\abs{y-x}\leq 30\sqrt{d}L_n$, $$P_{y,\omega}(X^*_{L_n^2}\geq v)\leq \exp(-\frac{v}{D_n}).$$\end{con}

The following proposition establishes, on the event $A_n$, the desired comparison between the continuous exit time $\tau^\epsilon$ and discrete stopping time $\tau^{\epsilon,n}_1$ with respect to $P_{x,\omega}$ for large $n$ and and on scales $\epsilon$ satisfying $L_n\leq \frac{1}{\epsilon}<L_{n+1}$.

\begin{prop}\label{end_time}  Assume (\ref{steady}) and (\ref{constants}).  For each $n\geq \overline{m}$ sufficiently large, for every $\epsilon>0$ satisfying $L_n\leq\frac{1}{\epsilon}<L_{n+1}$ and for every $\omega\in A_n$, for $C>0$ independent of $n$, $$\sup_{x\in \overline{U}/\epsilon}P_{x,\omega}(\tau^\epsilon-\tau^{\epsilon,n}_1\geq L_{n-1}^2)\leq CL_{n-1}^{-\zeta}.$$\end{prop}

\begin{proof}  Fix $n_1\geq 0$ as in Proposition \ref{end_Brownian} such that, for each $n\geq n_1$, $$2\tilde{D}_{n-\overline{m}}\leq \frac{r_0L_n}{2}.$$ This condition guarantees that Proposition \ref{end_Brownian} is satisfied for every $n\geq n_1$.  Furthermore, fix $n_2\geq 0$ such that, whenever $n\geq n_2$, $$L_{n+1}\overline{U}\subset[-\frac{1}{2}L_{n+2}^2, \frac{1}{2}L_{n+2}^2]^d,$$ which guarantees, whenever $n\geq n_2$ and $L_n\leq \frac{1}{\epsilon}<L_{n+2}$, the containment $\overline{U}/\epsilon\subset[-\frac{1}{2}L_{n+2}^2, \frac{1}{2}L_{n+2}^2]^d$ and therefore, for every $x\in\overline{U}/\epsilon$, the conclusion of Corollary \ref{couple_cor}.

Henceforth, fix $n\geq \max(n_1, n_2, \overline{m})$, $\epsilon>0$ satisfying $L_n\leq\frac{1}{\epsilon}<L_{n+1}$, $\omega\in A_n$ and $x\in \overline{U}/\epsilon$.  Recall the measure $Q_{n,x}$ defining the Markov chain $(X_k,\overline{X}_k)$ on $(\mathbb{R}^d\times\mathbb{R}^d)^{\mathbb{N}}$, and which effectively acts in its respective coordinates as a discrete version of the process in random environment or Brownian motion with variance $\alpha_{n-\overline{m}}$ along the sequence $\left\{kL_{n-\overline{m}}^2\right\}_{k\geq 0}$ in time.  Let $C_n$ denotes the event $$C_n=(\abs{X_k-\overline{X}_k}\geq L_{n-\overline{m}}\;|\;\textrm{for some}\;0\leq k\leq 2(\frac{L_{n+2}}{L_{n-\overline{m}}})^2\;),$$ and recall, owing to Corollary \ref{couple_cor}, for $C>0$ independent of $n$, \begin{equation}\label{end_time_1}Q_{x,n}(C_n)\leq C\tilde{\kappa}_{n-\overline{m}}L_{n-\overline{m}}^{16a-\delta}.\end{equation}

Let $\tilde{\tau}^\epsilon$ denotes the discrete $\C([0,\infty);\mathbb{R}^d)$ stopping time \begin{equation}\label{end_time_2} \tilde{\tau}^\epsilon=\inf\left\{\;kL_{n-\overline{m}}^2\geq0\;|\;X_{kL_{n-\overline{m}}^2}\notin U/\epsilon\;\right\}.\end{equation}  It follows by definition that $\tau^\epsilon\leq \tilde{\tau}^\epsilon$.  Furthermore, the definition of $Q_{n,y}$ and the Markov property imply that \begin{multline}\label{end_time_6}P_{x,\omega}(\tilde{\tau}^\epsilon-\tau^{\epsilon,n}_1\geq L_{n-1}^2)=Q_{n,x}(\tilde{T}^\epsilon-T^{\epsilon,n}_1\geq (\frac{L_{n-1}}{L_{n-\overline{m}}})^2) \\ =Q_{n,x}(\tilde{T}^\epsilon-T^{\epsilon,n}_1\geq (\frac{L_{n-1}}{L_{n-\overline{m}}})^2,C_n)+Q_{n,x}(\tilde{T}^\epsilon-T^{\epsilon,n}_1\geq (\frac{L_{n-1}}{L_{n-\overline{m}}})^2,C_n^c),\end{multline} where, in the final two terms of the equality, the stopping times are defined by $$T^{\epsilon,n}_1=\inf\left\{\;k\geq 0\;|\;(X_k,\overline{X}_k)\;\textrm{satisfies}\;d(X_k,(U/\epsilon)^c)\leq \tilde{D}_{n-\overline{m}}\;\right\},$$ which is merely the analogue of $\tau^{\epsilon,n}_1$ defined for the first coordinate of $(X_k, \overline{X}_k)$, and $$\tilde{T}^\epsilon=\inf\left\{\;k\geq 0\;|\;(X_k,\overline{X}_k)\;\textrm{satisfies}\;X_k\notin U\;\right\},$$ which is the analogue of $\tilde{\tau}^\epsilon$ defined for the first coordinate of $(X_k,\overline{X}_k)$.

In view of (\ref{end_time_1}), the first term of (\ref{end_time_6}) is bounded, for $C>0$ independent of $n$, by \begin{equation}\label{end_time_3} Q_{n,x}(\tilde{T}^\epsilon-T^{\epsilon,n}_1\geq (\frac{L_{n-1}}{L_{n-\overline{m}}})^2, C_n)\leq Q_{n,x}(C_n)\leq C\tilde{\kappa}_{n-\overline{m}}L_n^{16a-\delta}.\end{equation}  The second event is decomposed one step further as \begin{multline}\label{end_time_4} Q_{n,x}(\tilde{T}^\epsilon-T^{\epsilon,n}_1\geq (\frac{L_{n-1}}{L_{n-\overline{m}}})^2, C_n^c)=Q_{n,x}(\tilde{T}^\epsilon-T^{\epsilon,n}_1\geq (\frac{L_{n-1}}{L_{n-\overline{m}}})^2, C_n^c,T^{\epsilon,n}_1> (\frac{L_{n-1}}{L_{n-\overline{m}}})^2) \\ +Q_{n,x}(\tilde{T}^\epsilon-T^{\epsilon,n}_1\geq (\frac{L_{n-1}}{L_{n-\overline{m}}})^2, C_n^c,T^{\epsilon,n}_1\leq (\frac{L_{n-1}}{L_{n-\overline{m}}})^2).\end{multline}  Proposition \ref{exit_main}, and in particular line (\ref{exit_main_5}) which applies equally to the discrete sequence, since $L_{n-\overline{m}}^2$ divides $L_{n+2}^2$ according to the choice (\ref{L}), implies that the first term of (\ref{end_time_4}) can be estimated, for $C>0$ independent of $n$, by \begin{equation}\label{end_time_7}Q_{n,x}(\tilde{T}^\epsilon-T^{\epsilon,n}_1\geq (\frac{L_{n-1}}{L_{n-\overline{m}}})^2,  C_n^c, T^{\epsilon,n}_1> (\frac{L_{n-1}}{L_{n-\overline{m}}})^2)\leq Q_{n,x}(T^{\epsilon,n}_1>(\frac{L_{n-1}}{L_{n-\overline{m}}})^2)\leq CL_n^{-da}.\end{equation}  It remains to bound $$Q_{n,x}(\tilde{T}^\epsilon-T^{\epsilon,n}_1\geq (\frac{L_{n-1}}{L_{n-\overline{m}}})^2, C_n^c,T^{\epsilon,n}_1\leq (\frac{L_{n-1}}{L_{n-\overline{m}}})^2).$$

Define the discrete stopping time $$\overline{T}^{\epsilon,n}_2=\inf\left\{\;k\geq 0\;|\;(X_k,\overline{X}_k)\;\;\textrm{satisfies}\;\;\d(\overline{X}_k,(U/\epsilon))\geq \tilde{D}_{n-\overline{m}}\right\},$$ which is the analogue of $\tau^{\epsilon,n}_2$ defined for the second coordinate of the process $(X_k,\overline{X}_k)$ on $(\mathbb{R}^d\times\mathbb{R}^d)^{\mathbb{N}}$.  On the event $C_n^c$, for every $0\leq k\leq 2(\frac{L_{n+2}}{L_{n-\overline{m}}})^2$, whenever $$\d(X_k,(U/\epsilon)^c)\leq\tilde{D}_{n-\overline{m}},\;\;\textrm{it follows that}\;\;\d(\overline{X}_k,(U/\epsilon)^c)\leq \tilde{D}_{n-\overline{m}}+L_{n-\overline{m}}\leq 2\tilde{D}_{n-\overline{m}},$$ and whenever $$\d(\overline{X}_k,(U/\epsilon))\geq \tilde{D}_{n-\overline{m}},\;\;\textrm{it follows that}\;\;\d(X_k,(U/\epsilon))\geq \tilde{D}_{n-\overline{m}}-L_{n-\overline{m}}> 0.$$  Therefore, on the event $(C_n^c, T^{\epsilon,n}_1\leq (\frac{L_{n-1}}{L_{n-\overline{m}}})^2)$, since it follows from the definitions that $$\d(X_{T^{\epsilon,n}_1},(U/\epsilon)^c)\leq\tilde{D}_{n-\overline{m}}\;\;\textrm{and}\;\;\d(\overline{X}_{\overline{T}^{\epsilon,n}_2},(U/\epsilon))\geq\tilde{D}_{n-\overline{m}},$$ the Markov property, the definition of $Q_{n,x}$ and Proposition \ref{end_Brownian} imply, for $C>0$ independent of $n$, \begin{multline}\label{end_time_5} Q_{n,x}(\tilde{T}^\epsilon-T^{\epsilon,n}_1\geq (\frac{L_{n-1}}{L_{n-\overline{m}}})^2, C_n^c, T^{\epsilon,n}_1\leq (\frac{L_{n-1}}{L_{n-\overline{m}}})^2)\leq \sup_{\d(x,(U/\epsilon)^c)\leq 2\tilde{D}_{n-\overline{m}}}Q_{n,x}(\overline{T}^{\epsilon,n}_2\geq (\frac{L_{n-1}}{L_{n-\overline{m}}})^2) \\ =\sup_{\d(x,(U/\epsilon)^c)\leq 2\tilde{D}_{n-\overline{m}}}W^{n-\overline{m}}_x(\tau_2^{\epsilon,n}\geq L_{n-1}^2)\leq CL_{n-1}^{-\zeta}.\end{multline}

Therefore, owing to the choice of $\zeta>0$ in (\ref{end_zeta}) and since $\tau^\epsilon\leq \tilde{\tau}^\epsilon$ from (\ref{end_time_2}), the string of inequalities (\ref{end_time_6}), (\ref{end_time_3}), (\ref{end_time_4}), (\ref{end_time_7}) and (\ref{end_time_5}) imply, for $C>0$ independent of $n$, $$P_{x,\omega}(\tau^\epsilon-\tau^{\epsilon,n}_1\geq L_{n-1}^2)\leq P_{x,\omega}(\tilde{\tau}^\epsilon-\tau^{\epsilon,n}_1\geq L_{n-1}^2)\leq CL_{n-1}^{-\zeta},$$ which, since $x\in\overline{U}/\epsilon$, $n$ sufficiently large, $L_n\leq\frac{1}{\epsilon}<L_{n+1}$ and $\omega\in A_n$ were arbitrary, completes the argument.\end{proof}

The subsets $A_n$ now come to define the event on which the conclusion of Theorem \ref{intro_main} is obtained.  Recall Proposition \ref{prob_probability}, which states that, for each $n\geq \overline{m}$, for $C\geq 0$ independent of $n$, $$\mathbb{P}(\Omega\setminus A_n)\leq CL_n^{2d(1+a)^2-\frac{1}{2}M_0},$$  and notice that the definition of $L_n$ in (\ref{L}) and the negative exponent $2d(1+a)^2-\frac{1}{2}M_0<0$ guarantee the sum $$\sum_{n=\overline{m}}^\infty\mathbb{P}(\Omega\setminus A_n)\leq C\sum_{n=\overline{m}}^\infty L_n^{2d(1+a)^2-\frac{1}{2}M_0}<\infty.$$  Therefore, using the Borel-Cantelli lemma, let $\Omega_0\subset\Omega$ denote the subset of full probability \begin{equation}\label{end_mainevent} \Omega_0=\left\{\;\omega\in\Omega\;|\;\textrm{There exists}\;\overline{n}=\overline{n}(\omega)\;\textrm{such that}\;\omega\in A_n\;\textrm{for all}\;n\geq \overline{n}.\;\right\}.\end{equation}  Observe here that the set $\Omega_0$ is independent of $U$ and the boundary data.

Before shortly proceeding with the proof, it is convenient to recall some notation.  For each $n\geq 0$, let $u_n$ denote the solution \begin{equation}\label{end_hom_approx}\left\{\begin{array}{ll} \frac{\alpha_n}{2}\Delta u_n=0 & \textrm{on}\;\;U, \\ u_n=f(x) & \textrm{on}\;\;\partial U,\end{array}\right.\end{equation} and let $\overline{u}$ denote the solution \begin{equation}\label{end_hom}\left\{\begin{array}{ll} \Delta\overline{u}=0 & \textrm{on}\;\;U, \\ \overline{u}=f(x) & \textrm{on}\;\;\partial U.\end{array}\right.\end{equation}  The following proposition follows immediately by uniqueness and Theorem \ref{effectivediffusivity}, and states simply that the exit distribution from $U$ of a Brownian motion is independent of its (non-vanishing) variance, which corresponds to a time-change.

\begin{prop}\label{end_obvious}  Assume (\ref{steady}) and (\ref{constants}).  For each $n\geq 0$, $$u_n=\overline{u}\;\;\textrm{on}\;\;\overline{U}.$$\end{prop}

Similarly, for each $\epsilon>0$ and $\omega\in\Omega$, let $u^\epsilon$ denote the solution \begin{equation}\label{end_eq}\left\{\begin{array}{ll} \frac{1}{2}\tr(A(\frac{x}{\epsilon},\omega)D^2u^\epsilon)+\frac{1}{\epsilon}b(\frac{x}{\epsilon},\omega)\cdot Du^\epsilon=0 & \textrm{on}\;\;U, \\ u^\epsilon=f(x) & \textrm{on}\;\;\partial U.\end{array}\right.\end{equation}  The following theorem proves that, on the event $\Omega_0$, as $\epsilon\rightarrow 0$ the solutions $u^\epsilon$ converge uniformly to $\overline{u}$ on $\overline{U}$ whenever the boundary data is the restriction of a smooth function defined on the whole space.  Namely, \begin{equation}\label{end_compact} \textrm{assume}\;\;f\in C^\infty_c(\mathbb{R}^d).\end{equation}  This restriction will be removed by a standard approximation argument in the section's final theorem.

\begin{thm}\label{end} Assume (\ref{steady}), (\ref{constants}) and (\ref{end_compact}).  For every $\omega\in\Omega_0$, the solutions of (\ref{end_hom}) and (\ref{end_eq}) satisfy$$\lim_{\epsilon\rightarrow 0}\norm{u^\epsilon-\overline{u}}_{L^\infty(\overline{U})}=0.$$\end{thm}

\begin{proof}  Fix $\omega\in\Omega_0$ and $n_1\geq \overline{m}$ such that $\omega\in A_n$ for every $n\geq n_1$ and such that, whenever $n\geq n_1$, $$L_{n+1}U\subset[-\frac{1}{2}L_{n+2}^2, \frac{1}{2}L_{n+2}^2]^d$$ and the conditions of Propositions \ref{end_Brownian} and \ref{end_time} are satisfied.  Furthermore, choose $\epsilon_0\geq 0$ such that, whenever $0<\epsilon<\epsilon_0$, it follows that $L_n\leq \frac{1}{\epsilon}<L_{n+1}$ implies $n\geq n_1$.

The proof will rely upon the previously encountered continuous and discrete $\C([0,\infty);\mathbb{R}^d)$ stopping times, defined for each $\epsilon>0$ and $n\geq 0$, $$\tau^\epsilon=\inf\left\{\;t\geq 0\;|\;\epsilon X_t\notin U\right\}\;\;\textrm{and}\;\;\tau^{\epsilon,n}_1=\inf\left\{\;kL_{n-\overline{m}}^2\geq 0\;|\;\d(X_{kL_{n-\overline{m}}^2},(U/\epsilon)^c)\leq\tilde{D}_{n-\overline{m}}\;\right\},$$ and will use the representation $$u^\epsilon(x)=E_{\frac{x}{\epsilon},\omega}(f(\epsilon X_{\tau^\epsilon}))\;\;\textrm{on}\;\;\overline{U},$$ where $E_{\frac{x}{\epsilon},\omega}$ denotes the expectation corresponding to the $P_{\frac{x}{\epsilon},\omega}$ describing diffusion associated to the unscaled generator (\ref{end_gen}).

\emph{The Discrete Approximation:}  Fix $0<\epsilon<\epsilon_0$, the $n\geq 0$ satisfying $L_n\leq\frac{1}{\epsilon}<L_{n+1}$ and $x\in \overline{U}$.  First, decompose the representation in terms of the discrete approximation by \begin{equation}\label{end_1}u^\epsilon(x)=E_{\frac{x}{\epsilon},\omega}(f(\epsilon X_{\tau^\epsilon}))=E_{\frac{x}{\epsilon},\omega}(f(\epsilon X_{\tau^\epsilon})-f(\epsilon X_{\tau^{\epsilon,n}_1}))+E_{\frac{x}{\epsilon},\omega}(f(\epsilon X_{\tau^{\epsilon,n}_1})).\end{equation}  It will be shown that the first term of (\ref{end_1}) is negligible.

Decompose the expectation of the difference like \begin{multline}\label{end_01}E_{\frac{x}{\epsilon},\omega}(f(\epsilon X_{\tau^\epsilon})-f(\epsilon X_{\tau^{\epsilon,n}_1})) = E_{\frac{x}{\epsilon},\omega}(f(\epsilon X_{\tau^\epsilon})-f(\epsilon X_{\tau^{\epsilon,n}_1})), \tau^\epsilon+L_{n-\overline{m}}^2>\tau^{\epsilon,n}_1) \\ +E_{\frac{x}{\epsilon},\omega}(f(\epsilon X_{\tau^\epsilon})-f(\epsilon X_{\tau^{\epsilon,n}_1})),\tau^\epsilon+L_{n-\overline{m}}^2\leq \tau^{\epsilon,n}_1).\end{multline}  The event $\tau^\epsilon+L_{n-\overline{m}}^2$ implies by definition that the process beginning from $X_{\tau^\epsilon}$ travels further than $\tilde{D}_{n-\overline{m}}$ in time $L_{n-\overline{m}}^2$.  Therefore, the Markov property, $\omega\in A_n$, the choice of $n_1$ and the exponential estimated provided by Control \ref{end_localization} imply that \begin{multline}\label{end_2}\abs{E_{\frac{x}{\epsilon},\omega}(f(\epsilon X_{\tau^\epsilon})-f(\epsilon X_{\tau^{\epsilon,n}_1}),\tau^\epsilon+L_{n-\overline{m}}^2\leq \tau^{\epsilon,n}_1)}\leq \\ 2\norm{f}_{L^\infty(\mathbb{R}^d)}\exp(-\frac{\tilde{D}_{n-\overline{m}}}{D_{n-\overline{m}}})\leq \norm{f}_{L^\infty(\mathbb{R}^d)}\exp(-\kappa_{n-\overline{m}}).\end{multline}  The first term (\ref{end_01}) is further decomposed in the form \begin{multline}\label{end_4} E_{\frac{x}{\epsilon},\omega}(f(\epsilon X_{\tau^\epsilon})-f(\epsilon X_{\tau^{\epsilon,n}_1}),\tau^\epsilon+L_{n-\overline{m}}^2> \tau^{\epsilon,n}_1)= \\ E_{\frac{x}{\epsilon},\omega}(f(\epsilon X_{\tau^\epsilon})-f(\epsilon X_{\tau^{\epsilon,n}_1}), -L_{n-\overline{m}}^2\leq \tau^\epsilon-\tau^{\epsilon_n}_1< L_{n-1}^2) \\ +E_{\frac{x}{\epsilon},\omega}(f(\epsilon X_{\tau^\epsilon})-f(\epsilon X_{\tau^{\epsilon,n}_1}), \tau^\epsilon-\tau^{\epsilon,n}_1\geq L_{n-1}^2).\end{multline}  In view of Proposition \ref{end_time} and the choice of $0<\epsilon<\epsilon_0$, the second term (\ref{end_4}) is bounded, for $C>0$ independent of $n$, by \begin{multline}\label{end_400} \abs{E_{\frac{x}{\epsilon},\omega}(f(\epsilon X_{\tau^\epsilon})-f(\epsilon X_{\tau^{\epsilon,n}_1}), \tau^\epsilon-\tau^{\epsilon,n}_1\geq L_{n-1}^2)}\leq \\ 2\norm{f}_{L^\infty(\mathbb{R}^d)}P_{\frac{x}{\epsilon},\omega}(\tau^\epsilon-\tau^{\epsilon,n}_1\geq L_{n-1}^2)\leq C\norm{f}_{L^\infty(\mathbb{R}^d)}L_{n-1}^{-\zeta}.\end{multline}  The first term of (\ref{end_4}) is separated into the events \begin{multline}\label{end_6}E_{\frac{x}{\epsilon},\omega}(f(\epsilon X_{\tau^\epsilon})-f(\epsilon X_{\tau^{\epsilon,n}_1}), -L_{n-\overline{m}}^2\leq \tau^\epsilon-\tau^{\epsilon_n}_1< L_{n-1}^2)= \\ E_{\frac{x}{\epsilon},\omega}(f(\epsilon X_{\tau^\epsilon})-f(\epsilon X_{\tau^{\epsilon,n}_1}), 0< \tau^\epsilon-\tau^{\epsilon_n}_1< L_{n-1}^2) \\ +E_{\frac{x}{\epsilon},\omega}(f(\epsilon X_{\tau^\epsilon})-f(\epsilon X_{\tau^{\epsilon,n}_1}), 0\leq \tau^{\epsilon,n}_1-\tau^\epsilon< L_{n-\overline{m}}^2).\end{multline}  The Markov property and $f\in C^{\infty}_c(\mathbb{R}^d)$ imply \begin{multline*}  \abs{E_{\frac{x}{\epsilon},\omega}(f(\epsilon X_{\tau^\epsilon})-f(\epsilon X_{\tau^{\epsilon,n}_1}),  0< \tau^\epsilon-\tau^{\epsilon_n}_1< L_{n-1}^2)}\leq \\ \norm{Df}_{L^\infty(\mathbb{R}^d)}\epsilon \tilde{D}_{n-1}E_{\frac{x}{\epsilon},\omega}(P_{X_{\tau^{\epsilon,n}_1,\omega}}(X^*_{L_{n-1}^2}\leq \tilde{D}_{n-1})) \\ +2\norm{f}_{L^\infty(\mathbb{R}^d)}E_{\frac{x}{\epsilon},\omega}(P_{X_{\tau^{\epsilon,n}_1,\omega}}(X^*_{L_{n-1}^2}> \tilde{D}_{n-1})),\end{multline*} which, since $\omega\in A_n$ and $\frac{1}{L_{n+1}}<\epsilon\leq\frac{1}{L_n},$ it follows from Control \ref{end_localization} that, for $C>0$ independent of $n$, \begin{multline}\label{end_7}\abs{E_{\frac{x}{\epsilon},\omega}(f(\epsilon X_{\tau^\epsilon})-f(\epsilon X_{\tau^{\epsilon,n}_1}), 0< \tau^\epsilon-\tau^{\epsilon_n}_1< L_{n-1}^2)} \leq \\ C\norm{Df}_{L^\infty(\mathbb{R}^d)}\frac{\tilde{D}_{n-1}}{L_n}+C\norm{f}_{L^\infty(\mathbb{R}^d)}\exp(-\kappa_{n-1}).\end{multline}  And, the identical argument at scale $L_{n-\overline{m}}$ implies that the second term of (\ref{end_6}) satisfies, for $C>0$ independent of $n$, \begin{multline}\label{end_8} \abs{E_{\frac{x}{\epsilon},\omega}(f(\epsilon X_{\tau^\epsilon})-f(\epsilon X_{\tau^{\epsilon,n}_1}), 0\leq \tau^{\epsilon,n}_1-\tau^\epsilon< L_{n-\overline{m}}^2)}\leq \\ C\norm{Df}_{L^\infty(\mathbb{R}^d)}\frac{\tilde{D}_{n-\overline{m}}}{L_n}+C\norm{f}_{L^\infty(\mathbb{R}^d)}\exp(-\kappa_{n-\overline{m}}).\end{multline}

Therefore, inequalities (\ref{end_7}) and (\ref{end_8}) bound (\ref{end_6}) and show, for $C>0$ independent of $n$, \begin{multline*} \abs{E_{\frac{x}{\epsilon},\omega}(f(\epsilon X_{\tau^\epsilon})-f(\epsilon X_{\tau^{\epsilon,n}_1}), -L_{n-\overline{m}}^2\leq \tau^\epsilon-\tau^{\epsilon_n}_1< L_{n-1}^2)}\leq  \\ C\norm{Df}_{L^\infty(\mathbb{R}^d)}\frac{\tilde{D}_{n-1}}{L_n}+C\norm{f}_{L^\infty(\mathbb{R}^d)}\exp(-\kappa_{n-\overline{m}}).\end{multline*}  Then, combining this inequality with (\ref{end_400}) to bound (\ref{end_4}), for $C>0$ independent of $n$, since there exists $C>0$ such that $\exp(-\kappa_{n-\overline{m}})\leq CL_{n-1}^{-\zeta}$ for every $n\geq\overline{m}$,\begin{equation}\abs{E_{\frac{x}{\epsilon},\omega}(f(\epsilon X_{\tau^\epsilon})-f(\epsilon X_{\tau^{\epsilon,n}_1}),\tau^\epsilon+L_{n-\overline{m}}^2> \tau^{\epsilon,n}_1)}\leq C\norm{f}_{L^\infty(\mathbb{R}^d)}L_{n-1}^{-\zeta}+C\norm{Df}_{L^\infty(\mathbb{R}^d)}\frac{\tilde{D}_{n-1}}{L_n}.\end{equation}  And, using this inequality with (\ref{end_2}), the expectation of the difference (\ref{end_01}) can be estimated in the form, for $C>0$ independent of $n$, again using the fact that there exists $C>0$ independent of $n$ such that $\exp(-\kappa_{n-\overline{m}})\leq CL_{n-1}^{-\zeta}$ for all $n\geq\overline{m}$,\begin{equation}\label{end_9}\abs{E_{\frac{x}{\epsilon},\omega}(f(\epsilon X_{\tau^\epsilon})-f(\epsilon X_{\tau^{\epsilon,n}_1}))}\leq C\norm{f}_{L^\infty(\mathbb{R}^d)}L_{n-1}^{-\zeta}+C\norm{Df}_{L^\infty(\mathbb{R}^d)}\frac{\tilde{D}_{n-1}}{L_n}.\end{equation}

Therefore, in view of the decomposition (\ref{end_1}) and the estimate (\ref{end_9}), for $C>0$ independent of $n$, \begin{equation}\label{end_10} \abs{u^\epsilon(x)-E_{\frac{x}{\epsilon},\omega}(f(\epsilon X_{\tau^{\epsilon,n}_1})}\leq C\norm{f}_{L^\infty(\mathbb{R}^d)}L_{n-1}^{-\zeta}+C\norm{Df}_{L^\infty(\mathbb{R}^d)}\frac{\tilde{D}_{n-1}}{L_n}.\end{equation}  This estimate proves the efficacy of the discrete approximation defined by the stopping time $\tau_1^{\epsilon,n}$.  It will now be shown that the discretely stopped process is a good approximation of Brownian motion via the coupling estimate obtained in Corollay \ref{couple_cor}.

\emph{The Coupling:}  Recall the measure $Q_{n,\frac{x}{\epsilon}}$ defining the Markov chain $(X_k,\overline{X}_k)$ on $(\mathbb{R}^d\times\mathbb{R}^d)^{\mathbb{N}}$, and which effectively acts in the respective coordinates as discrete versions of the process in random environment and Brownian motion with variance $\alpha_{n-\overline{m}}$ along the sequence $\left\{kL_{n-\overline{m}}^2\right\}_{k\geq 0}$ in time.  Let $C_n$ denotes the event \begin{equation}\label{end_19}C_n=(\abs{X_k-\overline{X}_k}\geq L_{n-\overline{m}}\;|\;\textrm{for some}\;0\leq k\leq 2(\frac{L_{n+2}}{L_{n-\overline{m}}})^2\;),\end{equation} which, owing to Corollary \ref{couple_cor} and $\omega\in A_n$ with $n\geq n_1$, satisfies, for $C>0$ independent of $n$, \begin{equation}\label{end_20}Q_{n,\frac{x}{\epsilon}}(C_n)\leq C\tilde{\kappa}_{n-\overline{m}}L_{n-\overline{m}}^{16a-\delta}.\end{equation}  Furthermore, define as before the discrete stopping time $$T^{\epsilon,n}_1=\inf\left\{\;k\geq 0\;|\;(X_k,\overline{X}_k)\;\;\textrm{satisfies}\;d(X_k,(U/\epsilon)^c)\leq\tilde{D}_{n-\overline{m}}\;\right\},$$ which is the analogue of $\tau^{\epsilon,n}_1$ in the first coordinate.

The definition of $Q_{n,\frac{x}{\epsilon}}$ and the Markov property imply that, writing $E^{Q_{n,\frac{x}{\epsilon}}}$ for the expectation with respect to $Q_{n,\frac{x}{\epsilon}}$,\begin{equation}\label{end_11}E_{\frac{x}{\epsilon},\omega}(f(\epsilon X_{\tau^{\epsilon,n}_1})=E^{Q_{n,\frac{x}{\epsilon}}}(f(\epsilon X_{T^{\epsilon,n}_1}))=E^{Q_{n,\frac{x}{\epsilon}}}(f(\epsilon X_{T^{\epsilon,n}_1})-f(\epsilon \overline{X}_{T^{\epsilon,n}_1})))+E^{Q_{n,\frac{x}{\epsilon}}}(f(\epsilon \overline{X}_{T^{\epsilon,n}_1})).\end{equation}  As before, it will be shown that the expectation of the difference is negligible.  Decompose it in terms of the event $C_n$ to obtain \begin{multline}\label{end_21}E^{Q_{n,\frac{x}{\epsilon}}}(f(\epsilon X_{T^{\epsilon,n}_1})-f(\epsilon \overline{X}_{T^{\epsilon,n}_1})))=E^{Q_{n,\frac{x}{\epsilon}}}(f(\epsilon X_{T^{\epsilon,n}_1})-f(\epsilon \overline{X}_{T^{\epsilon,n}_1})), C_n) \\ +E^{Q_{n,\frac{x}{\epsilon}}}(f(\epsilon X_{T^{\epsilon,n}_1})-f(\epsilon \overline{X}_{T^{\epsilon,n}_1})), C_n^c).\end{multline}  The first term of (\ref{end_21}) is bounded using (\ref{end_20}) which implies, for $C>0$ independent of $n$, \begin{equation}\label{end_22} \abs{E^{Q_{n,\frac{x}{\epsilon}}}(f(\epsilon X_{T^{\epsilon,n}_1})-f(\epsilon \overline{X}_{T^{\epsilon,n}_1})), C_n)}\leq 2\norm{f}_{L^\infty(\mathbb{R}^d)}Q_{n,\frac{x}{\epsilon}}(C_n)\leq C\norm{f}_{L^\infty(\mathbb{R}^d)}\tilde{\kappa}_{n-\overline{m}}L_{n-\overline{m}}^{16a-\delta}.\end{equation}

The second term of (\ref{end_21}) is further decomposed in the form \begin{multline}\label{end_23} E^{Q_{n,\frac{x}{\epsilon}}}(f(\epsilon X_{T^{\epsilon,n}_1})-f(\epsilon \overline{X}_{T^{\epsilon,n}_1})), C_n^c)=E^{Q_{n,\frac{x}{\epsilon}}}(f(\epsilon X_{T^{\epsilon,n}_1})-f(\epsilon \overline{X}_{T^{\epsilon,n}_1})), C_n^c, T^{\epsilon,n}_1<(\frac{L_{n+2}}{L_{n-\overline{m}}})^2) \\ +E^{Q_{n,\frac{x}{\epsilon}}}(f(\epsilon X_{T^{\epsilon,n}_1})-f(\epsilon \overline{X}_{T^{\epsilon,n}_1})), C_n^c, T^{\epsilon,n}_1\geq (\frac{L_{n+2}}{L_{n-\overline{m}}})^2),\end{multline} observing here that $T^{\epsilon_n}_1=k$ corresponds to the process at time $kL_{n-\overline{m}}^2$.  The first term of (\ref{end_23}) is bounded using the definition of the set $C_n^c$, $T^{\epsilon,n}_1<(\frac{L_{n+2}}{L_{n-\overline{m}}})^2$ and $\epsilon\leq\frac{1}{L_n}$, which imply, for $C>0$ independent of $n$, \begin{multline}\label{end_24}\abs{E^{Q_{n,\frac{x}{\epsilon}}}(f(\epsilon X_{T^{\epsilon,n}_1})-f(\epsilon \overline{X}_{T^{\epsilon,n}_1})), C_n^c, T^{\epsilon,n}_1<(\frac{L_{n+2}}{L_{n-\overline{m}}})^2)}\leq \\ C\norm{Df}_{L^\infty(\mathbb{R}^d)}\epsilon L_{n-\overline{m}}\leq C\norm{Df}_{L^\infty(\mathbb{R}^d)}\frac{L_{n-\overline{m}}}{L_n}.\end{multline}  The second term of (\ref{end_23}) is bounded using the control for the exit time obtained in Proposition \ref{exit_main}, and in particular line (\ref{exit_main_5}) which applies equally to the discrete sequence since $L_{n-\overline{m}}^2$ divides $L_{n+2}^2$, to yield, for $C>0$ independent of $n$ and $\zeta>0$ defined in (\ref{end_zeta}), \begin{multline}\label{end_50}\abs{E^{Q_{n,\frac{x}{\epsilon}}}(f(\epsilon X_{T^{\epsilon,n}_1})-f(\epsilon \overline{X}_{T^{\epsilon,n}_1})), C_n^c, T^{\epsilon,n}_1\geq (\frac{L_{n+2}}{L_{n-\overline{m}}})^2)}\leq \\  2\norm{f}_{L^\infty(\mathbb{R}^d)}Q_{\frac{x}{\epsilon},n}(T^{\epsilon,n}_1\geq(\frac{L_{n+2}}{L_{n-\overline{m}}})^2)\leq C\norm{f}_{L^\infty(\mathbb{R}^d)}L_n^{-da}\leq C\norm{f}_{L^\infty(\mathbb{R}^d)}L_{n-1}^{-\zeta}.\end{multline}

Therefore, inequalities (\ref{end_24}) and (\ref{end_50}) bound (\ref{end_23}), for $C>0$ independent of $n$, by $$\abs{E^{Q_{n,\frac{x}{\epsilon}}}(f(\epsilon X_{T^{\epsilon,n}_1})-f(\epsilon \overline{X}_{T^{\epsilon,n}_1})), C_n^c)}\leq C\norm{f}_{L^\infty(\mathbb{R}^d)}L_{n-1}^{-\zeta}+C\norm{Df}_{L^\infty(\mathbb{R}^d)}\frac{L_{n-\overline{m}}}{L_n},$$ and together with the choice of $\zeta>0$ in (\ref{end_zeta}) and (\ref{end_22}), the expectation of the difference in (\ref{end_21}) can be estimated in the form, for $C>0$ independent of $n$, $$\abs{E^{Q_{n,\frac{x}{\epsilon}}}(f(\epsilon X_{T^{\epsilon,n}_1})-f(\epsilon \overline{X}_{T^{\epsilon,n}_1})))}\leq C\norm{f}_{L^\infty(\mathbb{R}^d)}L_{n-1}^{-\zeta}+C\norm{Df}_{L^\infty(\mathbb{R}^d)}\frac{L_{n-\overline{m}}}{L_n}.$$  And therefore, using (\ref{end_11}), for $C>0$ independent of $n$, \begin{equation}\label{end_25} \abs{E_{\frac{x}{\epsilon},\omega}(f(\epsilon X_{\tau^{\epsilon,n}_1})-E^{Q_{n,\frac{x}{\epsilon}}}(f(\epsilon \overline{X}_{T^{\epsilon,n}_1}))}\leq C\norm{f}_{L^\infty(\mathbb{R}^d)}L_{n-1}^{-\zeta}+C\norm{Df}_{L^\infty(\mathbb{R}^d)}\frac{L_{n-\overline{m}}}{L_n}.\end{equation}  It remains to recover the exit distribution of Brownian motion from the second term in the difference.

\emph{Recovering the Exit Distribution of Brownian Motion:}  The arguments here are essentially the unwinding, in terms of Brownian motion, of what led from (\ref{end_1}) to (\ref{end_10}).  Define the discrete stopping time $$\overline{T}^{\epsilon,n}_2=\inf\left\{\;k\geq 0\;|\;(X_k,\overline{X}_k)\;\;\textrm{satisfies}\;\d(\overline{X}_k,(U/\epsilon))\geq \tilde{D}_{n-\overline{m}}\;\right\},$$ which is the analogue of $\tau^{\epsilon,n}_2$ defined in (\ref{end_tau2}) for the second coordinate.  After performing decompositions analogous to (\ref{end_21}) and (\ref{end_23}), it follows by an identical argument that, for $C>0$ independent of $n$, \begin{equation}\label{end_26} \abs{E^{Q_{n,\frac{x}{\epsilon}}}(f(\epsilon \overline{X}_{T^{\epsilon,n}_1}))-E^{Q_{n,\frac{x}{\epsilon}}}(f(\epsilon \overline{X}_{T^{\epsilon,n}_1}), C_n^c, T^{\epsilon,n}_1<(\frac{L_{n+2}}{L_{n-\overline{m}}})^2)}\leq C\norm{f}_{L^\infty(\mathbb{R}^d)}L_{n-1}^{-\zeta}.\end{equation}  Notice that, since $T^{\epsilon,n}_1<(\frac{L_{n+2}}{L_{n-\overline{m}}})^2$, on the event $C_n^c$ \begin{equation}\label{end_27}\d(X_{T^{\epsilon,n}_1}, (U/\epsilon)^c)\leq\tilde{D}_{n-\overline{m}}\;\;\textrm{implies}\;\;\d(\overline{X}_{T^{\epsilon,n}_1}, (U/\epsilon)^c)\leq L_{n-\overline{m}}+\tilde{D}_{n-\overline{m}}\leq 2\tilde{D}_{n-\overline{m}}.\end{equation}  The second term of the difference in (\ref{end_26}) is then decomposed with respect to the stopping time $\overline{T}^{\epsilon,n}_2$ as in (\ref{end_1}), and takes the form \begin{multline}\label{end_29} E^{Q_{n,\frac{x}{\epsilon}}}(f(\epsilon \overline{X}_{T^{\epsilon,n}_1}), C_n^c, T^{\epsilon,n}_1<(\frac{L_{n+2}}{L_{n-\overline{m}}})^2)=E^{Q_{n,\frac{x}{\epsilon}}}(f(\epsilon \overline{X}_{T^{\epsilon,n}_1})-f(\epsilon \overline{X}_{\overline{T}^{\epsilon,n}_2}), C_n^c, T^{\epsilon,n}_1<(\frac{L_{n+2}}{L_{n-\overline{m}}})^2) \\ +E^{Q_{n,\frac{x}{\epsilon}}}(f(\epsilon \overline{X}_{\overline{T}^{\epsilon,n}_2}), C_n^c, T^{\epsilon,n}_1<(\frac{L_{n+2}}{L_{n-\overline{m}}})^2).\end{multline}  As before, the expectation of the difference is shown to be negligible.

Since, on the event $(C_n^c, (\frac{L_{n+2}}{L_{n-\overline{m}}})^2)$ it is necessarily the case that $T^{\epsilon,n}_1\leq \overline{T}^{\epsilon,n}_2$, the first term of (\ref{end_29}) is written\begin{multline}\label{end_28}E^{Q_{n,\frac{x}{\epsilon}}}(f(\epsilon \overline{X}_{T^{\epsilon,n}_1})-f(\epsilon \overline{X}_{\overline{T}^{\epsilon,n}_2}), C_n^c, T^{\epsilon,n}_1<(\frac{L_{n+2}}{L_{n-\overline{m}}})^2)= \\ E^{Q_{n,\frac{x}{\epsilon}}}(f(\epsilon \overline{X}_{T^{\epsilon,n}_1})-f(\epsilon \overline{X}_{\overline{T}^{\epsilon,n}_2}), C_n^c, T^{\epsilon,n}_1<(\frac{L_{n+2}}{L_{n-\overline{m}}})^2, \overline{T}^{\epsilon,n}_2-T^{\epsilon,n}_1>(\frac{L_{n-1}}{L_{n-\overline{m}}})^2) \\ +E^{Q_{n,\frac{x}{\epsilon}}}(f(\epsilon \overline{X}_{T^{\epsilon,n}_1})-f(\epsilon \overline{X}_{\overline{T}^{\epsilon,n}_2}), C_n^c, T^{\epsilon,n}_1<(\frac{L_{n+2}}{L_{n-\overline{m}}})^2, \overline{T}^{\epsilon,n}_2-T^{\epsilon,n}_1\leq (\frac{L_{n-1}}{L_{n-\overline{m}}})^2).\end{multline}  In view of (\ref{end_27}), the Markov property, the definition of $Q_{n,\frac{x}{\epsilon}}$ and Proposition \ref{end_Brownian}, the first term of (\ref{end_28}) is bounded, for $C>0$ independent of $n$, by \begin{multline}\label{end_30} \abs{E^{Q_{n,\frac{x}{\epsilon}}}(f(\epsilon \overline{X}_{T^{\epsilon,n}_1})-f(\epsilon \overline{X}_{\overline{T}^{\epsilon,n}_2}), C_n^c, T^{\epsilon,n}_1<(\frac{L_{n+2}}{L_{n-\overline{m}}})^2, \overline{T}^{\epsilon,n}_2-T^{\epsilon,n}_1>(\frac{L_{n-1}}{L_{n-\overline{m}}})^2)}\leq \\ 2\norm{f}_{L^\infty(\mathbb{R}^d)}\sup_{\d(x,(U/\epsilon)^c)\leq 2\tilde{D}_{n-\overline{m}}}W^{n-\overline{m}}_x(\tau_2^{\epsilon,n}\geq L_{n-1}^2)\leq C\norm{f}_{L^\infty(\mathbb{R}^d)}L_{n-1}^{-\zeta}.\end{multline}  Exponential estimates for Brownian motion will be used on scale $\tilde{D}_{n-1}$ to bound the second term of (\ref{end_28}), see \cite[Chapter~2, Proposition~1.8]{RY}, after performing the further decomposition \begin{multline}\label{end_31} E^{Q_{n,\frac{x}{\epsilon}}}(f(\epsilon \overline{X}_{T^{\epsilon,n}_1})-f(\epsilon \overline{X}_{\overline{T}^{\epsilon,n}_2}), C_n^c, T^{\epsilon,n}_1<(\frac{L_{n+2}}{L_{n-\overline{m}}})^2, \overline{T}^{\epsilon,n}_2-T^{\epsilon,n}_1\leq (\frac{L_{n+2}}{L_{n-\overline{m}}})^2)= \\ E^{Q_{n,\frac{x}{\epsilon}}}(f(\epsilon \overline{X}_{T^{\epsilon,n}_1})-f(\epsilon \overline{X}_{\overline{T}^{\epsilon,n}_2}), C_n^c, T^{\epsilon,n}_1<(\frac{L_{n+2}}{L_{n-\overline{m}}})^2, \\ \hspace{75mm} \overline{T}^{\epsilon,n}_2-T^{\epsilon,n}_1\leq (\frac{L_{n-1}}{L_{n-\overline{m}}})^2, \abs{\overline{X}_{T^{\epsilon,n}_1}- \overline{X}_{\overline{T}^{\epsilon,n}_2}}>\tilde{D}_{n-1}) \\ +E^{Q_{n,\frac{x}{\epsilon}}}(f(\epsilon \overline{X}_{T^{\epsilon,n}_1})-f(\epsilon \overline{X}_{\overline{T}^{\epsilon,n}_2}), C_n^c, T^{\epsilon,n}_1<(\frac{L_{n+2}}{L_{n-\overline{m}}})^2, \\ \overline{T}^{\epsilon,n}_2-T^{\epsilon,n}_1\leq (\frac{L_{n-1}}{L_{n-\overline{m}}})^2, \abs{\overline{X}_{T^{\epsilon,n}_1}- \overline{X}_{\overline{T}^{\epsilon,n}_2}}\leq\tilde{D}_{n-1}).\end{multline}  In analogy with (\ref{end_7}), the first term is bounded using the Markov property and exponential estimates for Brownian motion and the control of $\alpha_n$ provided by Theorem \ref{effectivediffusivity}, and the second term using the continuity of $f$ and $\epsilon\leq \frac{1}{L_n}$, to yield for $C,c>0$ independent of $n$, \begin{multline}\label{end_32}\abs{E^{Q_{n,\frac{x}{\epsilon}}}(f(\epsilon \overline{X}_{T^{\epsilon,n}_1})-f(\epsilon \overline{X}_{\overline{T}^{\epsilon,n}_2}), C_n^c, T^{\epsilon,n}_1<(\frac{L_{n+2}}{L_{n-\overline{m}}})^2, \overline{T}^{\epsilon,n}_2-T^{\epsilon,n}_1\leq (\frac{L_{n-1}}{L_{n-\overline{m}}})^2)}\leq  \\ C\norm{Df}_{L^\infty(\mathbb{R}^d)}\frac{\tilde{D}_{n-1}}{L_n}+C\norm{f}_{L^\infty(\mathbb{R}^d)}\exp(-c\kappa_{n-1}^2).\end{multline}  Therefore, combining (\ref{end_28}), (\ref{end_30}) and (\ref{end_32}), and using the fact that there exists $C>0$ independent of $n\geq 1$ such that $\exp(-c\kappa_{n-1}^2)\leq CL_{n-1}^{-\zeta}$, equation (\ref{end_29}) yields the estimate, for $C>0$ independent of $n$, \begin{multline}\label{end_33} \abs{E^{Q_{n,\frac{x}{\epsilon}}}(f(\epsilon \overline{X}_{T^{\epsilon,n}_1}), C_n^c, T^{\epsilon,n}_1<(\frac{L_{n+2}}{L_{n-\overline{m}}})^2)-E^{Q_{n,\frac{x}{\epsilon}}}(f(\epsilon \overline{X}_{\overline{T}^{\epsilon,n}_2}), C_n^c, T^{\epsilon,n}_1<(\frac{L_{n+2}}{L_{n-\overline{m}}})^2)}\leq \\ C\norm{Df}_{L^\infty(\mathbb{R}^d)}\frac{\tilde{D}_{n-1}}{L_n}+C\norm{f}_{L^\infty(\mathbb{R}^d)}L_{n-1}^{-\zeta}.\end{multline}  And, after repeating exactly the argument leading to (\ref{end_26}), for $C>0$ independent of $n$, \begin{equation}\label{end_34}\abs{E^{Q_{n,\frac{x}{\epsilon}}}(f(\epsilon \overline{X}_{\overline{T}^{\epsilon,n}_2}))-E^{Q_{n,\frac{x}{\epsilon}}}(f(\epsilon \overline{X}_{\overline{T}^{\epsilon,n}_2}), C_n^c, T^{\epsilon,n}_1<(\frac{L_{n+2}}{L_{n-\overline{m}}})^2)}\leq C\norm{f}_{L^\infty(\mathbb{R}^d)}L_{n-1}^{-\zeta}.\end{equation}  Then, combining (\ref{end_33}) and (\ref{end_34}) with (\ref{end_26}) yields, for $C>0$ independent of $n$, \begin{equation}\label{end_35} \abs{E^{Q_{n,\frac{x}{\epsilon}}}(f(\epsilon \overline{X}_{T^{\epsilon,n}_1}))-E^{Q_{n,\frac{x}{\epsilon}}}(f(\epsilon \overline{X}_{\overline{T}^{\epsilon,n}_2}))}\leq C\norm{Df}_{L^\infty(\mathbb{R}^d)}\frac{\tilde{D}_{n-1}}{L_n}+C\norm{f}_{L^\infty(\mathbb{R}^d)}L_{n-1}^{-\zeta}.\end{equation}  Since it follows from Markov property and the definition of $Q_{n,\frac{x}{\epsilon}}$ that $$E^{Q_{n,\frac{x}{\epsilon}}}(f(\epsilon \overline{X}_{\overline{T}^{\epsilon,n}_2}))=E^{W^{n-\overline{m}}_{\frac{x}{\epsilon}}}(f(\epsilon X_{\tau^{\epsilon,n}_2})),$$ equations (\ref{end_25}) and (\ref{end_35}) produce the estimate, for $C>0$ independent of $n$, \begin{equation}\label{end_36} \abs{E_{\frac{x}{\epsilon},\omega}(f(\epsilon X_{\tau^{\epsilon,n}_1})-E^{W^{n-\overline{m}}_{\frac{x}{\epsilon}}}(f(\epsilon X_{\tau^{\epsilon,n}_2}))}\leq C\norm{Df}_{L^\infty(\mathbb{R}^d)}\frac{\tilde{D}_{n-1}}{L_n}+C\norm{f}_{L^\infty(\mathbb{R}^d)}L_{n-1}^{-\zeta}.\end{equation}

\emph{Conclusion:}  It remains only to estimate the difference $$\abs{E^{W^{n-\overline{m}}_{\frac{x}{\epsilon}}}(f(\epsilon X_{\tau^{\epsilon,n}_2}))-E^{W^{n-\overline{m}}_{\frac{x}{\epsilon}}}(f(\epsilon X_{\tau^\epsilon}))}=\abs{E^{W^{n-\overline{m}}_{\frac{x}{\epsilon}}}(f(\epsilon X_{\tau^{\epsilon,n}_2})-f(\epsilon X_{\tau^\epsilon}))}.$$  Since it follows by definition that $\tau^\epsilon<\tau^{\epsilon,n}_2$, form the decomposition \begin{multline}\label{end_37} E^{W^{n-\overline{m}}_{\frac{x}{\epsilon}}}(f(\epsilon X_{\tau^{\epsilon,n}_2})-f(\epsilon X_{\tau^\epsilon}))=E^{W^{n-\overline{m}}_{\frac{x}{\epsilon}}}(f(\epsilon X_{\tau^{\epsilon,n}_2})-f(\epsilon X_{\tau^\epsilon}), \tau^{\epsilon,n}_2-\tau^\epsilon\leq L_{n-1}^2) \\ +E^{W^{n-\overline{m}}_{\frac{x}{\epsilon}}}(f(\epsilon X_{\tau^{\epsilon,n}_2})-f(\epsilon X_{\tau^\epsilon}), \tau^{\epsilon,n}_2-\tau^\epsilon>L_{n-1}^2).\end{multline}  Exponential estimates for Brownian motion on scale $\tilde{D}_{n-1}$, see \cite[Chapter~2, Proposition~1.8]{RY}, and $\epsilon<\frac{1}{L_n}$ imply that, for $C, c>0$ independent of $n$, and in exact analogy with the bound obtained in (\ref{end_32}), the first term of (\ref{end_37}) is bounded by \begin{multline}\label{end_38} \abs{E^{W^{n-\overline{m}}_{\frac{x}{\epsilon}}}(f(\epsilon X_{\tau^{\epsilon,n}_2})-f(\epsilon X_{\tau^\epsilon})), \tau^{\epsilon,n}_2-\tau^\epsilon\leq L_{n-1}^2)}\leq \\ C\norm{Df}_{L^\infty(\mathbb{R}^d)}\frac{\tilde{D}_{n-1}}{L_n}+C\norm{f}_{L^\infty(\mathbb{R}^d)}\exp(-c\kappa_{n-1}^2).\end{multline}  The second term of (\ref{end_37}) is handled similarly to (\ref{end_6}) but in the reverse order.  Here, \begin{multline}\label{end_39} \abs{E^{W^{n-\overline{m}}_{\frac{x}{\epsilon}}}(f(\epsilon X_{\tau^{\epsilon,n}_2})-f(\epsilon X_{\tau^\epsilon}), \tau^{\epsilon,n}_2-\tau^\epsilon>L_{n-1}^2)}\leq \\ 2\norm{f}_{L^\infty(\mathbb{R}^d)}W^{n-\overline{m}}_{\frac{x}{\epsilon}}(\tau^{\epsilon,n}_2-\tau^{\epsilon,n}_1>L_{n-1}^2-L_{n-\overline{m}}^2, \tau^{\epsilon,n}_1\leq \tau^\epsilon+L_{n-\overline{m}}^2) \\ +2\norm{f}_{L^\infty(\mathbb{R}^d)}W^{n-\overline{m}}_{\frac{x}{\epsilon}}(\tau^\epsilon+L_{n-\overline{m}}^2<\tau^{\epsilon,n}_1).\end{multline}  Since, for all $n\geq \overline{m}$, there exists $c_0>0$ satisfying $$L_{n-1}^2-L_{n-\overline{m}}^2\geq c_0L_{n-1}^2,$$ the proof of Proposition \ref{end_Brownian} implies that, for a larger $C>0$ independent of $n$, \begin{equation}\label{end_40} \sup_{\d(x,(U/\epsilon)^c)\leq 2\tilde{D}_{n-\overline{m}}}W^{n-\overline{m}}_x(\tau_2^{\epsilon,n}\geq c_0L_{n-1}^2)\leq CL_{n-1}^{-\zeta}.\end{equation}  Therefore, the Markov property, exponential estimates for Brownian motion on scale $\tilde{D}_{n-\overline{m}}$, see \cite[Chapter~2, Proposition~1.8]{RY}, (\ref{end_38}), (\ref{end_39}) and (\ref{end_40}) combine to bound (\ref{end_37}), using the fact that there exists $C>0$ such that $\exp(-\kappa_{n-1}^2)\leq \exp(-\kappa_{n-\overline{m}}^2)\leq CL_{n-1}^{-\zeta}$ for each $n\geq\overline{m}$, for $C>0$ independent of $n$, by \begin{equation}\label{end_41} \abs{E^{W^{n-\overline{m}}_{\frac{x}{\epsilon}}}(f(\epsilon X_{\tau^{\epsilon,n}_2})-f(\epsilon X_{\tau^\epsilon}))}\leq C\norm{f}_{L^\infty(\mathbb{R}^d)}L_{n-1}^{-\zeta}+C\norm{Df}_{L^\infty(\mathbb{R}^d)}\frac{\tilde{D}_{n-1}}{L_n}.\end{equation}

Finally, in view of (\ref{end_10}), (\ref{end_36}), (\ref{end_41}) and the triangle inequality, by the definition of $u^\epsilon$, $u_{n-\overline{m}}$, $\overline{u}$ and Proposition \ref{end_obvious}, for $C>0$ independent of $n$, \begin{multline}\label{end_42}\abs{E_{\frac{x}{\epsilon},\omega}(f(\epsilon X_{\tau^\epsilon})-E^{W^{n-\overline{m}}_{\frac{x}{\epsilon}}}(f(\epsilon X_{\tau^\epsilon}))}=\abs{u^\epsilon(x)-u_{n-\overline{m}}(x)}=\abs{u^\epsilon(x)-\overline{u}(x)}\leq \\  C\norm{Df}_{L^\infty(\mathbb{R}^d)}\frac{\tilde{D}_{n-1}}{L_n}+C\norm{f}_{L^\infty(\mathbb{R}^d)}L_{n-1}^{-\zeta}.\end{multline}  This, since $\epsilon\rightarrow 0$ implies $n\rightarrow\infty$ and the choice of constants (\ref{L}), (\ref{kappa}) and (\ref{D}) imply that $$\lim_{n\rightarrow\infty}(C\norm{Df}_{L^\infty(\mathbb{R}^d)}\frac{\tilde{D}_{n-1}}{L_n}+C\norm{f}_{L^\infty(\mathbb{R}^d)}L_{n-1}^{-\zeta})=0,$$ and because $x\in \overline{U}$, $\omega\in\Omega_0$ and $0<\epsilon<\epsilon_0$ were arbitrary, completes the proof.\end{proof}

The final theorem of this section extends Theorem \ref{end} to boundary data \begin{equation}\label{end_boundary_con} f\in\C(\partial U).\end{equation}  The proof follows by a standard approximation argument.

\begin{thm}\label{end_corollary} Assume (\ref{steady}), (\ref{constants}) and (\ref{end_boundary_con}).  For every $\omega\in\Omega_0$, the solutions of (\ref{end_hom}) and (\ref{end_eq}) satisfy $$\lim_{\epsilon\rightarrow 0}\norm{u^\epsilon-\overline{u}}_{L^\infty(\overline{U})}=0.$$\end{thm}

\begin{proof}  The Tietze Extension Theorem, see for instance Armstrong \cite[Page~40, Theorem~2.15]{Armstrong}, asserts that there exists a compactly supported extension $$\tilde{f}\in\BUC(\mathbb{R}^d)\;\;\textrm{satisfying}\;\;\tilde{f}|_{\partial U}=f.$$  Then, for each $\delta>0$, by convolution construct an $\tilde{f}^\delta\in\C^\infty_c(\mathbb{R}^d)$ satisfying $$\norm{\tilde{f}^\delta-\tilde{f}}_{L^\infty(\mathbb{R}^d)}\leq \delta,$$ and let $u^{\epsilon,\delta}$ denote the solution $$\left\{\begin{array}{ll} \frac{1}{2}\tr(A(\frac{x}{\epsilon},\omega)D^2 u^{\epsilon,\delta})+\frac{1}{\epsilon}b(\frac{x}{\epsilon},\omega)\cdot Du^{\epsilon,\delta}=0 & \textrm{on}\;\;U, \\ u^{\epsilon,\delta}=\tilde{f}^\delta(x) & \textrm{on}\;\;\partial U.\end{array}\right.$$  Similarly, let $\overline{u}^\delta$ denote the solution $$\left\{\begin{array}{ll} \Delta\overline{u}^{\delta}=0 & \textrm{on}\;\;U, \\ \overline{u}^\delta=\tilde{f}^\delta(x) & \textrm{on}\;\;\partial U.\end{array}\right.$$

The comparison principle and the triangle inequality imply that, for each $\omega\in \Omega$, $\delta>0$ and $\epsilon>0$, \begin{multline*}\norm{u^\epsilon-\overline{u}}_{L^\infty(\overline{U})}\leq \norm{u^\epsilon-u^{\epsilon,\delta}}_{L^\infty(\overline{U})}+\norm{u^{\epsilon,\delta}-\overline{u}^\delta}_{L^\infty(\overline{U})}+\norm{\overline{u}^\delta-\overline{u}}_{L^\infty(\overline{U})}\leq \\ 2\delta+\norm{u^{\epsilon,\delta}-\overline{u}^\delta}_{L^\infty(\overline{U})}.\end{multline*}  And, therefore, since $\tilde{f}^\delta$ satisfies the assumptions of Theorem \ref{end}, for every $\delta>0$ and $\omega\in\Omega_0$, $$\limsup_{\epsilon\rightarrow 0}\norm{u^\epsilon-\overline{u}}_{L^\infty(\overline{U})}\leq 2\delta,$$ which, since $\delta>0$ is arbitrary, completes the argument.  \end{proof}

\section{The Quantitative Estimate}\label{section_rate}

In the paper's final section, a rate for the convergence appearing in Theorem \ref{end_corollary} is first established for boundary data which is the restriction of a bounded, continuous function on $\mathbb{R}^d$.  That is, \begin{equation}\label{rate_restriction}\textrm{assume}\;f\in\BUC(\mathbb{R}^d),\end{equation} and write $\sigma_f:[0,\infty)\rightarrow[0,\infty)$ for the modulus of continuity \begin{equation}\label{rate_sigma} \abs{f(x)-f(y)}\leq\sigma_f(\abs{x-y})\;\;\textrm{for all}\;\;x,y\in\mathbb{R}^d.\end{equation}

Notice that, in the case $U=B_r$, which allows for an explicit radial extension, or whenever the domain $U$ is smooth, see the Product Neighborhood Theorem in Milnor \cite[Page~46]{Milnor}, every continuous function $f\in\C(\partial U)$, which is necessarily uniformly continuous by compactness, admits a continuous extension $\tilde{f}\in\BUC(\mathbb{R}^d)$ satisfying, for a constant $C=C(U)$ depending only upon the domain, $$\sigma_{\tilde{f}}(s)\leq \sigma_f(Cs)\;\;\textrm{for all}\;\;s\in[0,\infty).$$  And therefore, for sufficiently smooth domains, assumption (\ref{rate_restriction}) is always satisfied up to a domain dependent factor.

\begin{thm}\label{rate}  Assume (\ref{steady}), (\ref{constants}) and (\ref{rate_restriction}).  There exists $C>0$ such that, for every $\omega\in\Omega_0$, for all $\epsilon>0$ sufficiently small depending on $\omega$, for $\zeta>0$ defined in (\ref{end_zeta}), the solutions of (\ref{end_hom}) and (\ref{end_eq}) satisfy $$\norm{u^\epsilon-\overline{u}}_{L^\infty(\overline{U})}\leq C\norm{f}_{L^\infty(\mathbb{R}^d)} \epsilon^{\frac{\zeta}{2(1+a)^2}}+C\sigma_f(\epsilon^{\frac{a}{2(1+a)^2}}).$$\end{thm}

\begin{proof}  Fix $\omega\in\Omega_0$.  The only observation is that, in every step of the proof of Theorem \ref{end} involving the continuity of $f$, the Lipschitz estimates can be replaced by estimates using the modulus $\sigma_f$.  And, therefore, since $\omega\in\Omega_0$ and in view of the final estimate (\ref{end_42}), whenever $\epsilon>0$ is sufficiently small and $n\geq 0$ satisfies $L_n\leq \frac{1}{\epsilon}<L_{n+1}$, for $C>0$ independent of $n$ and $\omega$, $$\norm{u^\epsilon-\overline{u}}_{L^\infty(\overline{U})}\leq C\sigma_f(\frac{\tilde{D}_{n-1}}{L_n})+C\norm{f}_{L^\infty(\mathbb{R}^d)}L_{n-1}^{-\zeta}.$$  Then, it follows from the definition of the constants (\ref{L}), (\ref{kappa}) and (\ref{D}) that, for all $n\geq 0$ sufficiently large and whenever $L_n\leq \frac{1}{\epsilon}<L_{n+1}$, $$\frac{\tilde{D}_{n-1}}{L_n}\leq L_{n+1}^{-\frac{a}{2(1+a)^2}}\leq \epsilon^{\frac{a}{2(1+a)^2}}\;\;\textrm{and}\;\; L_{n-1}^{-\zeta}\leq L_{n+1}^{-\frac{\zeta}{2(1+a)^2}}\leq \epsilon^{\frac{\zeta}{2(1+a)^2}},$$ which, since $\omega\in\Omega_0$ was arbitrary, completes the argument.  \end{proof}

The paper's final theorem extends Theorem \ref{rate} to general continuous boundary data provided the domain $U$ is smooth.  In this case, \begin{equation}\label{rate_smooth_domain}\textrm{assume}\;f\in\C(\partial U)\;\textrm{and that the domain}\;U\;\textrm{is smooth}.\end{equation}  The proof is an immediate consequence of Theorem \ref{rate} and the remark immediately preceding.

\begin{thm}\label{rate_corollary}  Assume (\ref{steady}), (\ref{constants}) and (\ref{rate_smooth_domain}).  There exists $C>0$ and $C_1=C_1(U)>0$ such that, for every $\omega\in\Omega_0$, for all $\epsilon>0$ sufficiently small depending on $\omega$, for $\zeta>0$ defined in (\ref{end_zeta}),  the solutions of (\ref{end_hom}) and (\ref{end_eq}) satisfy $$\norm{u^\epsilon-\overline{u}}_{L^\infty(\overline{U})}\leq C\norm{f}_{L^\infty(\partial U)} \epsilon^{\frac{\zeta}{2(1+a)^2}}+C\sigma_f(C_1 \epsilon^{\frac{a}{2(1+a)^2}}).$$\end{thm}

\bibliography{Dirichlet}
\bibliographystyle{plain}

\end{document}